# CHARACTERIZATION OF THE ROBUST ISOLATED CALMNESS FOR A CLASS OF CONIC PROGRAMMING PROBLEMS

CHAO DING*, DEFENG SUN†, AND LIWEI ZHANG‡

**Abstract.** This paper is devoted to studying the robust isolated calmness of the Karush-Kuhn-Tucker (KKT) solution mapping for a large class of interesting conic programming problems (including most commonly known ones arising from applications) at a locally optimal solution. Under the Robinson constraint qualification, we show that the KKT solution mapping is robustly isolated calm if and only if both the strict Robinson constraint qualification and the second order sufficient condition hold. This implies, among others, that at a locally optimal solution the second order sufficient condition is needed for the KKT solution mapping to have the Aubin property.

**1. Introduction.** Let $\mathcal{X}$ and $\mathcal{Y}$ be two finite dimensional real Euclidean spaces each equipped with an inner product $\langle \cdot, \cdot \rangle$ and its induced norm $\|\cdot\|$. Consider the following canonically perturbed optimization problem:

(1)
$$\begin{aligned} \min \quad & f(x) - \langle a, x \rangle \\ \text{s.t.} \quad & G(x) + b \in \mathcal{K}, \end{aligned}$$

where $f : \mathcal{X} \to \Re$ and $G : \mathcal{X} \to \mathcal{Y}$ are twice continuously differentiable functions, $\mathcal{K} \subset \mathcal{Y}$ is a nonempty closed convex set, and $(a, b) \in \mathcal{X} \times \mathcal{Y}$ is the perturbation parameter.

For each given $(a, b) \in \mathcal{X} \times \mathcal{Y}$, we use $X(a, b)$ to denote the set of all locally optimal solutions of problem (1). A point $x \in X(a, b)$ is said to be isolated if there exists an open neighborhood $\mathcal{V}$ of $x$ such that $X(a, b) \cap \mathcal{V} = \{x\}$. Let $\Phi(a, b)$ be the set of all feasible points of problem (1) with a given $(a, b)$, i.e.,

(2) $$\Phi(a, b) := \{x \in \mathcal{X} \mid G(x) + b \in \mathcal{K}\}, \quad (a, b) \in \mathcal{X} \times \mathcal{Y}.$$

Let $L : \mathcal{X} \times \mathcal{Y} \to \Re$ be the Lagrangian function of problem (1) defined by

(3) $$L(x; y) := f(x) + \langle y, G(x) \rangle, \quad (x, y) \in \mathcal{X} \times \mathcal{Y}.$$

For any $y \in \mathcal{Y}$, denote the derivative of $L(\cdot; y)$ at $x \in \mathcal{X}$ by $L'_x(x; y)$ and denote the adjoint of $L'_x(x; y)$ by $\nabla_x L(x; y)$. For a given perturbation parameter $(a, b)$, the Karush-Kuhn-Tucker (KKT) optimality condition for problem (1) takes the following form:

(4) $$\begin{cases} a = \nabla_x L(x; y), \\ b \in -G(x) + \partial \sigma(y, \mathcal{K}) \end{cases} \iff \begin{cases} a = \nabla_x L(x; y), \\ y \in \mathcal{N}_{\mathcal{K}}(G(x) + b), \end{cases}$$

*Institute of Applied Mathematics, Academy of Mathematics and Systems Science, Chinese Academy of Sciences, Beijing, P.R. China (dingchao@amss.ac.cn). The research of this author was supported by the National Natural Science Foundation of China under project No. 11671387.

†Department of Mathematics and Risk Management Institute, National University of Singapore, 10 Lower Kent Ridge Road, Singapore (matsundf@nus.edu.sg). The research of this author was supported in part by the Academic Research Fund (Grant No. R-146-000-207-112).

‡School of Mathematical Sciences, Dalian University of Technology, Dalian, P.R. China. (lwzhang@dlut.edu.cn). The research of this author was supported by the National Natural Science Foundation of China under project No. 91330206 and No. 11571059.

where $\sigma(y, \mathcal{K}) := \sup\{\langle y, z\rangle \mid z \in \mathcal{K}\}$ is the support function of $\mathcal{K}$, $\partial \sigma(y, \mathcal{K})$ is the sub-differential of $\sigma(\cdot, \mathcal{K})$ at $y$ and $\mathcal{N}_{\mathcal{K}}(z)$ is the normal cone of $\mathcal{K}$ at $z \in \mathcal{K}$ in the context of convex analysis [33]. For a given $(a, b) \in \mathcal{X} \times \mathcal{Y}$, the set of all solutions $(x, y)$ to the KKT system (4) is denoted by $S_{\text{KKT}}(a, b)$. We use $X_{\text{KKT}}(a, b)$ to denote the set of all stationary points of problem (1) with $(a, b)$, i.e.,

$$X_{\text{KKT}}(a, b) := \{x \in \mathcal{X} \mid \text{there exists } y \in \mathcal{Y} \text{ such that (4) holds at } (x, y)\}.$$

The set of Lagrange multipliers associated with $(x, a, b)$ is defined by

(5)  $$M(x, a, b) := \{y \in \mathcal{Y} \mid (x, y) \in S_{\text{KKT}}(a, b)\}.$$

For each given $(a, b) \in \mathcal{X} \times \mathcal{Y}$, we say that the Robinson constraint qualification (RCQ) for problem (1) holds at a feasible point $x \in \Phi(a, b)$ if

(6) $$G'(\bar{x})\mathcal{X} + \mathcal{T}_{\mathcal{K}}(G(\bar{x})) = \mathcal{Y},$$

where for any subset $\Omega$ in a finite dimensional Euclidean space $\mathcal{E}$, the tangent cone [34, Definition 6.1] to $\Omega$ at $z \in \Omega$ is defined by

$$\mathcal{T}_\Omega(z) := \left\{d \in \mathcal{E} \mid \exists z^k \to z \text{ with } z^k \in \Omega \text{ and } t^k \downarrow 0 \text{ such that } (z^k - z)/t^k \to d\right\}.$$

For each given $(a, b)$, it is well-known (cf. e.g., [7, Theorem 3.9 & Proposition 3.17]) that the RCQ holds at a locally optimal solution $x \in \mathcal{X}$ if and only if $M(x, a, b)$ is a nonempty, convex, and compact subset of $\mathcal{Y}$.

In this paper, we mainly focus on an important property in perturbation theory for problem (1) at a locally optimal solution: the KKT solution mapping $S_{\text{KKT}}$ is locally nonempty-valued and isolated calm (see Definition 2 in Section 2). This property will be referred as the robust isolated calmness in this paper. When the set $\mathcal{K}$ in problem (1) is polyhedral, the theory on the robust isolated calmness is fairly complete even under a more general perturbation framework, e.g., the functions $f$ and $G$ in problem (1) have the parametric forms: $f(x, c)$ and $G(x, c)$, where $c \in \mathcal{Z}$ is another parameter in a finite dimensional real Euclidean space $\mathcal{Z}$ [32, 10, 18]. One natural question that one may ask is to what extent one can generalize the results from the polyhedral case to the non-polyhedral case. For example, Dontchev and Rockafellar [10] show that for the parametric nonlinear programming problem at a locally optimal solution, $S_{\text{KKT}}$ is robustly isolated calm if and only if both the strict Mangasarian-Fromovitz constraint qualification (MFCQ) and the second order sufficient condition (SOSC) hold. Does this result still hold in the non-polyhedral setting? In particular for the case $\mathcal{K} = \mathcal{S}^n_+$, the cone of $n$ by $n$ symmetric and positive semi-definite matrices in $\mathcal{S}^n$ (PSD cone in short), the Euclidean space of $n$ by $n$ symmetric matrices.

Before answering the above question, let us first consider the following convex quadratic semi-definite programming (SDP) example constructed by Bonnans and Shapiro.

EXAMPLE 1. *Bonnans and Shapiro [7, Example 4.54]:*

(7) $$\begin{aligned} \min \quad & x_1 + x_1^2 + x_2^2 \\ \text{s.t.} \quad & \text{Diag}(x) + \varepsilon A \in \mathcal{S}^2_+, \end{aligned}$$

*where $x = (x_1, x_2) \in \Re^2$, $\text{Diag}(x)$ is the $2 \times 2$ diagonal matrix whose $i$-th diagonal element is $x_i$, $i = 1, 2$, $A$ is a non-diagonal matrix in $\mathcal{S}^2$, and $\varepsilon$ is a scalar parameter. It is clear that the Slater condition holds. When $\varepsilon = 0$, the optimization problem*

(7) has the unique optimal solution $\bar{x} = (0, 0)$ with the unique Lagrange multiplier $\overline{Y} = \begin{bmatrix} -1 & 0 \\ 0 & 0 \end{bmatrix}$. It is easy to see that for any given $\varepsilon \geq 0$, problem (7) has a unique optimal solution $X(\varepsilon) = (\bar{x}_1(\varepsilon), \bar{x}_2(\varepsilon))$ with $\bar{x}_2(\varepsilon)$ of order $\varepsilon^{2/3}$ as $\varepsilon \to 0$.

Example 1 demonstrates that the KKT solution mapping $S_{\text{KKT}}$ of a convex quadratic SDP problem with a strongly convex objective function and with the Slater condition being satisfied can still fail to be calm. This is completely different from the polyhedral set case [32], where Robinson shows that if $\mathcal{K}$ is a polyhedral cone, the KKT solution mapping $S_{\text{KKT}}$ always possesses the calmness property if both the MFCQ and one strong form of the SOSC hold. Therefore, cautions must be taken when extensions are made from the polyhedral to the non-polyhedral case. One may argue that the unperturbed problem (7) in Example 1 is not a natural SDP problem as the unperturbed problem reduces to a nonlinear programming problem. However, it is not difficult to construct a non-lifted SDP example of similar behaviors (see Example 3 in Section 4). This also indicates the following fact: under the general non-polyhedral set setting, the uniqueness of Lagrange multipliers and the SOSC even in the strong form do not imply the calmness, let alone the isolated calmness, of $S_{\text{KKT}}$. Of course, for nonlinear programming problems at a locally optimal solution, since the strict MFCQ is equivalent to the uniqueness of Lagrange multipliers [21, Proposition 1.1], the aforementioned result of Dontchev and Rockafellar says that $S_{\text{KKT}}$ is robustly isolated calm if and only if the uniqueness of Lagrange multipliers and the SOSC hold.

The literature on perturbation analysis of general optimization problems is enormous, and even a short summary about the most important results achieved would be far beyond our reach. For recent works covering many topics in perturbation analysis, one may refer to [6, 7, 11, 13, 19, 25, 34] and references therein. Here, we only touch those results that are mostly relevant to the research conducted in this paper. When the non-polyhedral set $\mathcal{K}$ is the second-order cone (SOC) or the PSD cone, the characterizations of the strong regularity of the KKT point (see Robinson [31] for the definition) are established in [5] and [36], respectively. In [5], Bonnans and Ramírez show that for a locally optimal solution of the nonlinear SOC programming problem, under the RCQ, the KKT solution mapping is strongly regular if and only if the strong SOSC and the constraint non-degeneracy hold. In [36], Sun shows that for a locally optimal solution of the nonlinear SDP problem, the strong SOSC and the constraint non-degeneracy, the non-singularity of Clarke's Jacobian of the KKT system and the strong regularity are all equivalent under the RCQ. The strong regularity always implies the robust isolated calmness, but not vice versa. So, in order for characterizing the robust isolated calmness, we need weaker conditions than the strong SOSC and the constraint non-degeneracy.

From King and Rockafellar [17] and Levy [22], we know that the isolated calmness of $S_{\text{KKT}}$ for problem (1) with the non-polyhedral set $\mathcal{K}$ can be characterized by the non-singularity of its graphical derivative. Under the assumption that the constraint non-degeneracy (see (12) for the definition) holds at a stationary point, not necessarily a locally optimal solution, to problem (1), Mordukhovich et al. [26, 27] establish an explicit formula for the graphical derivative of $S_{\text{KKT}}$ and thus a characterization of the isolated calmness of $S_{\text{KKT}}$ in terms of the non-singularity of its graphical derivative. The constraint non-degeneracy assumption is crucial for the analysis conducted in [26, 27], in particular [26, Lemma 3.1] and [27, Proposition 3.1].

In general, the non-singularity of the graphical derivative of $S_{\text{KKT}}$ may be difficult to verify. However, recently there are some good progresses in addressing this nonsingularity issue if the concerned point is a locally optimal solution. In [37], Zhang and Zhang show that the strict Robinson constraint qualification (SRCQ) (see (11) for the definition) and the SOSC condition imply the non-singularity of the graphical derivative of $S_{\rm KKT}$ for the nonlinear SDP problem at a locally optimal solution. Han et al. [16] show that for the nonlinear SDP problem, the isolated calmness of $S_{\rm KKT}$ implies the SRCQ. Moreover, complete characterizations of the isolated calmness of $S_{\rm KKT}$ are provided for convex composite quadratic SDP problems [16, Theorem 5.1]. In [23], Liu and Pan extend several results in [16] and [37] to the matrix optimization problem constrained by the epigraph cone of the Ky Fan matrix $k$-norm. For the nonlinear SOC programming problem at a locally optimal solution, Zhang et al. [38] show that the isolated calmness of $S_{\rm KKT}$ is equivalent to the SRCQ and the SOSC.

In this paper, assuming that the set $\mathcal{K}$ in problem (1) belongs to the class of $C^2$-cone reducible sets (Definition 6), we study the continuity of $X_{\rm KKT}$ and the robust isolated calmness of $S_{\rm KKT}$. It is worth mentioning that the class of $C^2$-cone reducible sets is rich. It includes notably all the polyhedral convex sets and many non-polyhedral sets such as the SOC, the PSD cone [35, 7] and the epigraph cone of the Ky Fan matrix $k$-norm [8]. Moreover, the Cartesian product of $C^2$-cone reducible sets is also $C^2$-cone reducible (see, e.g., [7, 35]). Under the $C^2$-cone reducibility assumption, we will provide a sufficient condition for the stationary mapping $X_{\rm KKT}$ of the canonically perturbed problem (1) to be continuous and isolated. Furthermore, for the $C^2$-cone reducible case, under the RCQ, we present a complete characterization of the robust isolated calmness of the KKT solution mapping $S_{\rm KKT}$, namely, the robust isolated calmness of $S_{\rm KKT}$ holds at a locally optimal solution if and only if the corresponding SRCQ and the SOSC hold. Additionally, we construct an example (Example 4) to show that our conditions for characterizing the robust isolated calmness of $S_{\rm KKT}$ studied in this paper are strictly weaker than those for characterizing the Aubin property (see Definition 3 in Section 2) of the KKT system (4). Thus, by combining with [20, Theorem 1] established recently by Klatte and Kummer, we conclude that under the $C^2$-cone reducible assumption, at a locally optimal solution to problem (1), the constraint non-degeneracy and the SOSC are both necessary for $S_{\rm KKT}$ to have the Aubin property.

The remaining parts of this paper are organized as follows. In Section 2, we introduce some definitions and preliminary results on variational analysis. In Section 3, we study the continuity and the isolatedness of the stationary mapping $X_{\rm KKT}$ of problem (1). The characterization of the robust isolated calmness of the KKT solution mapping $S_{\rm KKT}$ and its implication for the Aubin property are provided in Section 4. We conclude our paper in Section 5.

**2. Preliminaries.** Firstly, let us recall some common notions and definitions related to set-valued mappings. Let $\mathcal{E}$ be a finite dimensional real Euclidean space. Let $\mathcal{E}$ and $\mathcal{F}$ be two finite dimensional real Euclidean spaces and $\Psi : \mathcal{E} \rightrightarrows \mathcal{F}$ be a set-valued mapping with $(\bar{p}, \bar{q}) \in {\rm gph}\,\Psi$, i.e., $\bar{q} \in \Psi(\bar{p})$, where ${\rm gph}\,\Psi$ denotes the graph of $\Psi$. Let $\mathbb{B}_\mathcal{F}$ be the unit ball in $\mathcal{F}$. The set-valued mapping $\Psi$ is said to be lower semi-continuous at $\bar{p}$ for $\bar{q}$ if for any open neighborhood $\mathcal{V}$ of $\bar{q}$ there exists an open neighborhood $\mathcal{U}$ of $\bar{p}$ such that

$$\emptyset \neq \Psi(p) \cap \mathcal{V} \quad \forall p \in \mathcal{U}.$$

The mapping $\Psi$ is said to be upper semi-continuous (in Berge's sense [2]) at $\bar{p}$ if for any open set $\mathcal{O} \supset \Psi(\bar{p})$ there exists an open neighborhood $\mathcal{U}$ such that for any $p \in \mathcal{U}$, $\Psi(p) \subset \mathcal{O}$. Furthermore, if $\Psi$ is lower semi-continuous at $(\bar{p}, \bar{q})$ and is upper semi-continuous at $\bar{p}$, then $\Psi$ is said to be continuous at $(\bar{p}, \bar{q}) \in \operatorname{gph} \Psi$.

For the set-valued mapping $\Psi$, we are interested in the following three Lipschitz-like properties: the calmness, isolated calmness and Aubin property.

DEFINITION 1. *The set-valued mapping $\Psi : \mathcal{E} \rightrightarrows \mathcal{F}$ is said to be calm at $\bar{p}$ if there exist a constant $\kappa > 0$ and an open neighborhood $\mathcal{U}$ of $\bar{p}$ such that*

$$\Psi(p) \subset \Psi(\bar{p}) + \kappa \|p - \bar{p}\| \mathbb{B}_{\mathcal{F}} \quad \forall\, p \in \mathcal{U}.$$

DEFINITION 2. *The set-valued mapping $\Psi : \mathcal{E} \rightrightarrows \mathcal{F}$ is said to be isolated calm at $\bar{p}$ for $\bar{q}$ if there exist a constant $\kappa > 0$ and open neighborhoods $\mathcal{U}$ of $\bar{p}$ and $\mathcal{V}$ of $\bar{q}$ such that*

$$(8) \qquad \Psi(p) \cap \mathcal{V} \subset \{\bar{q}\} + \kappa \|p - \bar{p}\| \mathbb{B}_{\mathcal{F}} \quad \forall\, p \in \mathcal{U}.$$

*Moreover, $\Psi$ is said to be robustly isolated calm*[1] *at $\bar{p}$ for $\bar{q}$ if (8) holds and for each $p \in \mathcal{U}$, $\Psi(p) \cap \mathcal{V} \neq \emptyset$.*

DEFINITION 3. *The set-valued mapping $\Psi : \mathcal{E} \rightrightarrows \mathcal{F}$ has the Aubin property at $\bar{p}$ for $\bar{q}$ if there exist a constant $\kappa > 0$ and open neighborhoods $\mathcal{U}$ of $\bar{p}$ and $\mathcal{V}$ of $\bar{q}$ such that*

$$\Psi(p) \cap \mathcal{V} \subset \Psi(p') + \kappa \|p - p'\| \mathbb{B}_{\mathcal{F}} \quad \forall\, p, p' \in \mathcal{U}.$$

The calmness for the set-valued mapping $\Psi$ given in Definition 1 comes from [34, 9(30)] and it was called "upper Lipschitzian" by Robinson [30]. The isolated calmness for the set-valued mapping $\Psi$ was called differently in the literature, e.g., the local upper Lipschitz continuity in [10, 22], to distinguish it from Robinson's definition of upper Lipschitz continuity [30]. The property defined by Definition 3 was designated "pseudo-Lipschitzian" by Aubin [1].

REMARK 1. *It is worth noting that the set-valued mapping $\Psi$ is isolated calm at $\bar{p}$ for $\bar{q}$ does not implies $\Psi$ is robust isolated calm at $\bar{p}$ for $\bar{q}$ in general (see, for instance, [27, Example 6.4]). However, if $\Psi$ is lower semi-continuous at $(\bar{p}, \bar{q})$, then it follows from the definition that the robust isolated calmness of $\Psi$ also holds.*

The graphical derivative [34, Definition 8.33] of $\Psi$ at $\bar{p}$ for $\bar{q}$ is a set-valued mapping $D\Psi(\bar{p}|\bar{q}) : \mathcal{E} \rightrightarrows \mathcal{F}$ defined by

$$D\Psi(\bar{p}|\bar{q})(u) := \{v \in \mathcal{F} \mid (u, v) \in \mathcal{T}_{\operatorname{gph} \Psi}(\bar{p}, \bar{q})\},$$

which is a convenient tool for investigating the isolated calmness property. In fact, we have the following basic characterization of the isolated calmness of the set-valued mapping $\Psi$ at $\bar{p}$ for $\bar{q}$.

LEMMA 4 (King and Rockafellar [17], Levy [22]). *Let $(\bar{p}, \bar{q}) \in \operatorname{gph} \Psi$. Then $\Psi$ is isolated calm at $\bar{p}$ for $\bar{q}$ if and only if $\{0\} = D\Psi(\bar{p}|\bar{q})(0)$.*

In order to study the relationship between the isolated calmness and the Aubin property of the set-valued mappings, we need the following result.

LEMMA 5 (Fusek [14]). *Suppose that $F : \mathcal{E} \to \mathcal{E}$ is locally Lipschitz continuous near $\bar{q} \in \mathcal{E}$ and that $F$ is directionally differentiable at $\bar{q}$. If $F^{-1}$ has the Aubin property at $\bar{p} := F(\bar{q})$ for $\bar{q}$, then there exists an open neighborhood $\mathcal{V}$ of $\bar{q}$ such that $F^{-1}(\bar{p}) \cap \mathcal{V} = \{\bar{q}\}$.*

---

[1] If $\Psi = \Theta^{-1}$ is the inverse of a given set-valued mapping $\Theta : \mathcal{Y} \rightrightarrows \mathcal{X}$, then the robust isolated calmness of $\Psi$ at $\bar{p}$ for $\bar{q}$ is equivalent to the upper regularity [19] of $\Theta$ at $\bar{q}$ for $\bar{p}$.

It then follows from Lemma 5 that for a function $F : \mathcal{E} \to \mathcal{E}$ satisfying the assumptions in this lemma, the Aubin property of $F^{-1}$ implies the isolated calmness of $F^{-1}$.

Suppose that $\bar{x}$ is a feasible solution to problem (1) with $(a, b) = (0, 0)$. The critical cone $\mathcal{C}(\bar{x})$ of (1) with $(a, b) = (0, 0)$ at $\bar{x}$ is defined by

$$\text{(9)} \qquad \mathcal{C}(\bar{x}) := \{d \in \mathcal{X} \mid G'(\bar{x})d \in \mathcal{T}_\mathcal{K}(G(\bar{x})),\ f'(\bar{x})d \leq 0\}.$$

If $\bar{x}$ is a stationary point of problem (1) with $(a, b) = (0, 0)$ and $\bar{y} \in M(\bar{x}, 0, 0)$, then

$$\begin{aligned}\mathcal{C}(\bar{x}) &= \{d \in \mathcal{X} \mid G'(\bar{x})d \in \mathcal{T}_\mathcal{K}(G(\bar{x})),\ f'(\bar{x})d = 0\} \\ &= \{d \in \mathcal{X} \mid G'(\bar{x})d \in \mathcal{C}_\mathcal{K}(G(\bar{x}), \bar{y})\},\end{aligned}$$

where for any $A \in \mathcal{K}$, $\mathcal{C}_\mathcal{K}(A, B)$ is the critical cone of $\mathcal{K}$ at $A$ with respect to $B \in \mathcal{N}_\mathcal{K}(A)$ defined by

$$\text{(10)} \qquad \mathcal{C}_\mathcal{K}(A, B) := \mathcal{T}_\mathcal{K}(A) \cap B^\perp,$$

and for any $s \in \mathcal{Y}$, $s^\perp := \{z \in \mathcal{Y} \mid \langle z, s \rangle = 0\}$. Since $\mathcal{K}$ is convex, it is clear that for each $\bar{y} \in M(\bar{x}, 0, 0)$, the critical cone $\mathcal{C}_\mathcal{K}(G(\bar{x}), \bar{y})$ is indeed a closed convex cone.

The SRCQ is said to hold for problem (1) with $(a, b) = (0, 0)$ at $\bar{x}$ with respect to $\bar{y} \in M(\bar{x}, 0, 0) \neq \emptyset$ if

$$\text{(11)} \qquad G'(\bar{x})\mathcal{X} + \mathcal{T}_\mathcal{K}(G(\bar{x})) \cap \bar{y}^\perp = \mathcal{Y}.$$

It follows from [7, Proposition 4.50] that the set of Lagrange multipliers $M(\bar{x}, 0, 0)$ is a singleton if the SRCQ holds. The constraint non-degeneracy, introduced by Robinson [31], is said to hold at $\bar{x}$ if

$$\text{(12)} \qquad G'(\bar{x})\mathcal{X} + \text{lin}\,(\mathcal{T}_\mathcal{K}(G(\bar{x}))) = \mathcal{Y},$$

where $\text{lin}\,(\mathcal{T}_\mathcal{K}(G(\bar{x})))$ is the lineality space of $\mathcal{T}_\mathcal{K}(G(\bar{x}))$, i.e., the largest linear space in $\mathcal{T}_\mathcal{K}(G(\bar{x}))$. It is well-known (cf. [7, Proposition 4.73]) that the constraint non-degeneracy is stronger than the SRCQ since $\text{lin}\,(\mathcal{T}_\mathcal{K}(G(\bar{x}))) \subset \mathcal{T}_\mathcal{K}(G(\bar{x})) \cap \bar{y}^\perp$.

Through the whole paper, we always assume that the set $\mathcal{K}$ has the following $C^2$-cone reducibility.

DEFINITION 6 ([7, Definition 3.135]). *The closed convex set $\mathcal{K}$ is said to be $C^2$-cone reducible at $\overline{A} \in \mathcal{K}$, if there exist a open neighborhood $\mathcal{W} \subset \mathcal{Y}$ of $\overline{A}$, a pointed closed convex cone $\mathcal{Q}$ (a cone is said to be pointed if and only if its lineality space is the origin) in a finite dimensional space $\mathcal{Z}$ and a twice continuously differentiable mapping $\Xi : \mathcal{W} \to \mathcal{Z}$ such that: (i) $\Xi(\overline{A}) = 0 \in \mathcal{Z}$; (ii) the derivative mapping $\Xi'(\overline{A}) : \mathcal{Y} \to \mathcal{Z}$ is onto; (iii) $\mathcal{K} \cap \mathcal{W} = \{A \in \mathcal{W} \mid \Xi(A) \in \mathcal{Q}\}$. We say that $\mathcal{K}$ is $C^2$-cone reducible if $\mathcal{K}$ is $C^2$-cone reducible at every $\overline{A} \in \mathcal{K}$.*

Let $Z$ be a closed set in $\mathcal{Y}$. Recall that the inner and outer second order tangent sets ([7, (3.49) and (3.50)]) to the given closed set $Z$ in the direction $h \in \mathcal{Y}$ can be defined, respectively, by

$$\mathcal{T}_Z^{i,2}(z, h) := \{w \in \mathcal{Y} \mid \text{dist}(z + th + \tfrac{1}{2}t^2 w, Z) = o(t^2),\ t \geq 0\}$$

and

$$\mathcal{T}_Z^2(z, h) := \{w \in \mathcal{Y} \mid \exists t_k \downarrow 0 \text{ such that } \text{dist}(z + t_k h + \tfrac{1}{2}t_k^2 w, Z) = o(t_k^2)\},$$

where for any $s \in \mathcal{Y}$, $\mathrm{dist}(s, Z) := \inf\{\|s - z\| \mid z \in Z\}$. Note that in general, $\mathcal{T}_Z^{i,2}(z, h) \neq \mathcal{T}_Z^2(z, h)$ even if $Z$ is convex ([7, Section 3.3]). However, it follows from [7, Proposition 3.136] that if $Z$ is a $C^2$-cone reducible convex set, then the equality always holds. In this case, $\mathcal{T}_Z^2(z, h)$ will be simply called the second order tangent set to $Z$ at $z \in Z$ in the direction $h \in \mathcal{Y}$.

For the $C^2$-cone reducible set $\mathcal{K}$, we have the following "no gap" second order conditions for problem (1) with $(a, b) = (0, 0)$ (see [4, Theorem 3.1 and Theorem 4.1]).

THEOREM 7. *Suppose that $\bar{x}$ is a locally optimal solution to problem* (1) *with $(a, b) = (0, 0)$ and the RCQ holds at $\bar{x}$. Then the following second-order necessary condition holds*

$$(13) \quad \sup_{y \in M(\bar{x},0,0)} \left\{ \langle d, \nabla_{xx}^2 L(\bar{x}; y) d \rangle - \sigma\left(y, \mathcal{T}_\mathcal{K}^2(G(\bar{x}), G'(\bar{x})d)\right) \right\} \geq 0 \quad \forall d \in \mathcal{C}(\bar{x}),$$

*where for any $y \in \mathcal{Y}$, $\nabla_{xx}^2 L(\bar{x}; y)$ is the Hessian of $L(\cdot; y)$ at $\bar{x}$. Conversely, suppose $\bar{x}$ is a stationary point of problem* (1) *with $(a, b) = (0, 0)$ and the RCQ holds at $\bar{x}$. Then the following second order sufficient condition (SOSC)*

$$(14) \quad \sup_{y \in M(\bar{x},0,0)} \left\{ \langle d, \nabla_{xx}^2 L(\bar{x}; y) d \rangle - \sigma\left(y, \mathcal{T}_\mathcal{K}^2(G(\bar{x}), G'(\bar{x})d)\right) \right\} > 0 \quad \forall d \in \mathcal{C}(\bar{x}) \setminus \{0\}$$

*is necessary and sufficient for the quadratic growth condition at the point $\bar{x}$ for problem* (1) *with respect to $(a, b) = (0, 0)$.*

Next, we list some results that are needed for our subsequent discussions from the standard reduction approach. For more details on the reduction approach, one may refer to [7, Section 3.4.4]. The following results on the representations of the normal cone and the "sigma term" of the $C^2$-cone reducible set $\mathcal{K}$ are stated in [7, (3.266) and (3.274)].

LEMMA 8. *Let $\overline{A} \in \mathcal{K}$ be given. Then, there exist an open neighborhood $\mathcal{W} \subset \mathcal{Y}$ of $\overline{A}$, a pointed closed convex cone $\mathcal{Q}$ in a finite dimensional space $\mathcal{Z}$ and a twice continuously differentiable function $\Xi : \mathcal{W} \to \mathcal{Z}$ satisfying conditions (i)-(iii) in Definition 6 such that for all $A \in \mathcal{W}$ sufficiently close to $\overline{A}$,*

$$(15) \quad \mathcal{N}_\mathcal{K}(A) = \Xi'(A)^* \mathcal{N}_\mathcal{Q}(\Xi(A)),$$

*where $\Xi'(A)^* : \mathcal{Z} \to \mathcal{Y}$ is the adjoint of $\Xi'(A)$. In particular, for any $\overline{B} \in \mathcal{N}_\mathcal{K}(\overline{A})$, there is a unique element $u$ in $\mathcal{N}_\mathcal{Q}(\Xi(\overline{A}))$ such that $\overline{B} = \Xi'(\overline{A})^* u$, denoted by $(\Xi'(\overline{A})^*)^{-1} \overline{B}$. Furthermore, we have for any $D \in \mathcal{C}_\mathcal{K}(\overline{A}, \overline{B})$,*

$$(16) \quad \sigma(\overline{B}, \mathcal{T}_\mathcal{K}^2(\overline{A}, D)) = -\langle (\Xi'(\overline{A})^*)^{-1} \overline{B}, \Xi''(\overline{A})(D, D) \rangle.$$

For a feasible solution $\bar{x}$ to problem (1) with $(a, b) = (0, 0)$, let $\mathcal{W}$, $\mathcal{Q}$ and $\Xi$ be the open neighborhood of $G(\bar{x})$, the pointed closed convex cone and the twice continuously differentiable function defined in Lemma 8, respectively, with respect to $G(\bar{x}) \in \mathcal{K}$. Since $G$ is continuous, we know that there exist open neighborhoods $\mathcal{U}$ of the origin in $\mathcal{X} \times \mathcal{Y}$ and $\mathcal{V}$ of $\bar{x}$ such that $G(x) + b \in \mathcal{W}$ for any $(x, a, b) \in \mathcal{V} \times \mathcal{U}$. Consequently, problem (1) is locally equivalent to the following reduced problem:

$$(17) \quad \begin{aligned} \min \quad & f(x) - \langle a, x \rangle \\ \text{s.t.} \quad & \mathcal{G}(x, b) \in \mathcal{Q}, \end{aligned}$$

where $\mathcal{G}(x, b) := \Xi(G(x) + b)$ for each $(x, a, b) \in \mathcal{V} \times \mathcal{U}$, in the sense that the sets of optimal solutions to (1) and (17) restricted to $\mathcal{V}$ are the same. Moreover, for $(a, b) = (0, 0)$, it is known from [7, Section 3.4.4] that the RCQ for problem (1) holds at the feasible point $\bar{x}$ if and only if the RCQ for problem (17) holds at $\bar{x}$.

Next, we shall present some useful results about the directional derivative of the metric projection operator over the $C^2$-cone reducible set $\mathcal{K}$. Suppose that $\overline{B} \in \mathcal{N}_\mathcal{K}(\overline{A})$. Let $C := \overline{A} + \overline{B}$. Then we have $\overline{A} = \Pi_\mathcal{K}(C)$, where $\Pi_\mathcal{K} : \mathcal{Y} \to \mathcal{Y}$ is the metric projection operator over $\mathcal{K}$, i.e., for any $C \in \mathcal{Y}$,

$$\Pi_\mathcal{K}(C) := \operatorname{argmin}\left\{\frac{1}{2}\|Y - C\|^2 \mid Y \in \mathcal{K}\right\}.$$

Since $\mathcal{K}$ is $C^2$-cone reducible, we know from [3, Theorem 7.2] that $\Pi_\mathcal{K}$ is directionally differentiable at $C$ and the directional derivative $\Pi'_\mathcal{K}(C; H)$ for any direction $H \in \mathcal{Y}$ is the unique optimal solution to the following strongly convex optimization problem:

$$(18) \qquad \min\left\{\|D - H\|^2 - \sigma(\overline{B}, \mathcal{T}^2_\mathcal{K}(\overline{A}, D)) \mid D \in \mathcal{C}_\mathcal{K}(\overline{A}, \overline{B})\right\},$$

where $\mathcal{C}_\mathcal{K}(\overline{A}, \overline{B})$ is the critical cone of $\mathcal{K}$ at $\overline{A}$ with respect to $\overline{B}$ defined by (10). It follows from (16) in Lemma 8 that there exists a self-adjoint linear operator $\mathcal{H} : \mathcal{Y} \to \mathcal{Y}$ such that

$$(19) \qquad \Upsilon(D) := \langle D, \mathcal{H}(D) \rangle = -\sigma(\overline{B}, \mathcal{T}^2_\mathcal{K}(\overline{A}, D)) \geq 0 \quad \forall D \in \mathcal{C}_\mathcal{K}(\overline{A}, \overline{B}),$$

which means that that $\mathcal{H}$ is co-positive on the cone $\mathcal{C}_\mathcal{K}(\overline{A}, \overline{B})$. However, this does not mean that $\mathcal{H}$ is positive semi-definite on the whole space $\mathcal{Y}$. To overcome this difficulty, let us define $h : \mathcal{Y} \to (-\infty, \infty]$ by

$$(20) \qquad h(D) := \Upsilon(D) + \delta_{\mathcal{C}_\mathcal{K}(\overline{A}, \overline{B})}(D), \quad D \in \mathcal{Y},$$

where $\delta_{\mathcal{C}_\mathcal{K}(\overline{A}, \overline{B})}(\cdot)$ is the indicator function of the critical cone $\mathcal{C}_\mathcal{K}(\overline{A}, \overline{B})$.

PROPOSITION 9. *The function $h : \mathcal{Y} \to (-\infty, \infty]$ defined by (20) is a closed proper convex function and the sub-differential $\partial h(D)$ of $h$ at any $D \in \mathcal{C}_\mathcal{K}(\overline{A}, \overline{B})$ is given by*

$$\partial h(D) = \nabla \Upsilon(D) + \mathcal{N}_{\mathcal{C}_\mathcal{K}(\overline{A}, \overline{B})}(D).$$

*Proof.* First, by using (19) and (20), we know that

$$h(D) = -\sigma(\overline{B}, \mathcal{T}^2_\mathcal{K}(\overline{A}, D)) + \delta_{\mathcal{C}_\mathcal{K}(\overline{A}, \overline{B})}(D) \quad \forall D \in \mathcal{Y}.$$

Since $\delta_{\mathcal{C}_\mathcal{K}(\overline{A}, \overline{B})}(D) = 0$ for any $D \in \mathcal{C}_\mathcal{K}(\overline{A}, \overline{B})$ and the function $-\sigma(\overline{B}, \mathcal{T}^2_\mathcal{K}(\overline{A}, \cdot))$ is a closed proper convex function [4, Lemma 4.1], it follows that $h(\cdot)$ is also a closed proper convex function on $\mathcal{Y}$. Moreover, we know from [34, Proposition 8.12] (or [25, Theorem 1.93]) that $\partial h(D) = \partial_L h(D)$ for any $D \in \mathcal{C}_\mathcal{K}(\overline{A}, \overline{B})$, where $\partial_L h(D)$ is the limiting sub-differential of $h$ at $D$ (cf. e.g., [25, Definition 1.77]). Thus, it follows from the sum rule [25, Proposition 1.107 (ii)] that

$$\partial h(D) = \partial_L h(D) = \nabla \Upsilon(D) + \mathcal{N}_{\mathcal{C}_\mathcal{K}(\overline{A}, \overline{B})}(D) \quad \forall D \in \mathcal{C}_\mathcal{K}(\overline{A}, \overline{B}).$$

This completes the proof of this proposition. □

For any given nonempty convex cone $K \subset \mathcal{Y}$, we use $K^\circ$ to denote the polar of $K$, i.e., $K^\circ := \{z \in \mathcal{Y} \mid \langle z, s \rangle \leq 0 \ \forall s \in K\}$. The following simple lemma is a generalization of [16, Lemmas 4.1 & 4.3].

LEMMA 10. *Let $C \in \mathcal{Y}$, $\overline{A} = \Pi_\mathcal{K}(C)$ and $\overline{B} = C - \overline{A}$.*

(i) Let $\triangle A, \triangle B \in \mathcal{Y}$. $\triangle A - \Pi'_{\mathcal{K}}(C; \triangle A + \triangle B) = 0$ if and only if

(21) $$\begin{cases} \triangle A \in \mathcal{C}_{\mathcal{K}}(\overline{A}, \overline{B}), \\ \triangle B - \frac{1}{2}\nabla\Upsilon(\triangle A) \in [\mathcal{C}_{\mathcal{K}}(\overline{A}, \overline{B})]^{\circ}, \\ \langle \triangle A, \triangle B \rangle = -\sigma(\overline{B}, \mathcal{T}^{2}_{\mathcal{K}}(\overline{A}, \triangle A)), \end{cases}$$

where $\Upsilon(\cdot)$ is the quadratic function defined by (19).

(ii) Let $\mathcal{A} : \mathcal{X} \to \mathcal{Y}$ be a linear operator. Then, the following two statements are equivalent:

(a) $\triangle B \in \mathcal{Y}$ is a solution to the following system of equations
$$\begin{cases} \mathcal{A}^{*}\triangle B = 0, \\ \Pi'_{\mathcal{K}}(C; \triangle B) = 0; \end{cases}$$

(b) $\triangle B \in \left[ \mathcal{A}\mathcal{X} + \mathcal{T}_{\mathcal{K}}(\overline{A}) \cap \overline{B}^{\perp} \right]^{\circ}$.

*Proof.* We first prove part (i). Let $h : \mathcal{Y} \to (-\infty, \infty]$ be defined by (20). From Proposition 9 we know that $h(\cdot)$ is a closed proper convex function and $\overline{D} \in \mathcal{C}_{\mathcal{K}}(\overline{A}, \overline{B})$ is the unique optimal solution to problem (18) with $H = \triangle A + \triangle B$ if and only if $\overline{D}$ is the unique optimal solution to the following strongly convex optimization problem

$$\min \left\{ \|D - (\triangle A + \triangle B)\|^{2} + h(D) \right\},$$

or equivalently,

$$0 \in 2(\overline{D} - (\triangle A + \triangle B)) + \nabla\Upsilon(\overline{D}) + \mathcal{N}_{\mathcal{C}_{\mathcal{K}}(\overline{A}, \overline{B})}(\overline{D}).$$

Therefore, by using the fact that $\mathcal{C}_{\mathcal{K}}(\overline{A}, \overline{B})$ is a closed convex cone, we can see that $\triangle A = \Pi'_{\mathcal{K}}(C; \triangle A + \triangle B)$ if and only if

$$\begin{cases} \triangle A \in \mathcal{C}_{\mathcal{K}}(\overline{A}, \overline{B}), \\ \triangle B - \frac{1}{2}\nabla\Upsilon(\triangle A) \in [\mathcal{C}_{\mathcal{K}}(\overline{A}, \overline{B})]^{\circ}, \\ \langle \triangle A, \triangle B - \frac{1}{2}\nabla\Upsilon(\triangle A) \rangle = 0, \end{cases}$$

which, together with the fact that for each $\triangle A \in \mathcal{C}_{\mathcal{K}}(\overline{A}, \overline{B})$,

$$\langle \triangle A, \nabla\Upsilon(\triangle A) \rangle = 2\Upsilon(\triangle A) = -2\sigma(B, \mathcal{T}^{2}_{\mathcal{K}}(A, \triangle A)),$$

show that (21) holds.

By noting $0 \in \mathcal{C}_{\mathcal{K}}(\overline{A}, \overline{B})$ and taking $\triangle A = 0$ in part (i), we obtain part (ii) immediately. $\square$

**3. The continuity and isolatedness of $X_{\mathrm{KKT}}$.** Let $x \in \mathcal{X}$ be a feasible solution to problem (1) with $(a, b) = (0, 0)$, i.e., $x \in \Phi(0, 0)$. It is well-known (cf. e.g., [7, Theorem 2.87] and [34, Theorem 9.43]) that the feasible point mapping $\Phi$ defined by (2) has the Aubin property at $(0, 0)$ for $x$ if the RCQ (6) holds at $x$. Therefore, the following result on the existence of local minimizers of problem (1) extends [10, Lemma 2.5] for nonlinear programming.

LEMMA 11. *Suppose that $\bar{x}$ is an isolated locally optimal solution of problem (1) with $(a, b) = (0, 0)$. If the RCQ (6) holds at $\bar{x}$, then the locally optimal solution mapping $X$ is lower semi-continuous at $(0, 0, \bar{x}) \in \mathrm{gph}\, X$.*

*Proof.* Since the feasible solution mapping $\Phi$ has the Aubin property at $(0,0)$ for $\bar{x}$, we know that there exist $\varepsilon_1 > 0$, $\varepsilon_2 > 0$ and $\kappa > 0$ such that for any $\|(a,b)\| < \varepsilon_1$, $\|(a',b')\| < \varepsilon_1$,

$$\Phi(a,b) \cap \{x \in \mathcal{X} \mid \|x - \bar{x}\| < \varepsilon_2\} \subset \Phi(a',b') + \kappa \|(a,b) - (a',b')\| \mathbb{B}_{\mathcal{Y}}. \tag{22}$$

Let $\mathcal{V}$ be an arbitrary open neighborhood of $\bar{x}$. Since $\bar{x}$ is an isolated locally optimal solution of problem (1) with $(a,b) = (0,0)$, we are able to choose $\tau \in (0, \varepsilon_2)$ sufficiently small such that $\{x \in \mathcal{X} \mid \|x - \bar{x}\| < \tau\} \subset \mathcal{V}$ and $\bar{x}$ is the unique minimizer of (1) with $(a,b) = (0,0)$ in $\{x \in \mathcal{X} \mid \|x - \bar{x}\| < \tau\}$. For the fixed $\tau$ and for any $\|(a,b)\| < \varepsilon_1$, define the following set-valued mapping

$$\Phi_\tau(a,b) := \Phi(a,b) \cap \{x \in \mathcal{X} \mid \|x - \bar{x}\| \le \tau + \kappa\|(a,b)\|\}.$$

Note that the graph of $\Phi_\tau$ is closed, which means that $\Phi_\tau$ is closed. Also, for any $\|(a,b)\| < \varepsilon_1$, $\Phi_\tau(a,b)$ belongs the compact subset $\{x \in \mathcal{X} \mid \|x - \bar{x}\| \le \varepsilon_2 + \kappa\varepsilon_1\}$. Therefore, we know (cf. e.g., [7, Lemma 4.3 (i)]) that $\Phi_\tau$ is upper semi-continuous at any $\|(a,b)\| < \varepsilon_1$.

On the other hand, for any $x \in \Phi_\tau(0,0) = \Phi(0,0) \cap \{x \in \mathcal{X} \mid \|x - \bar{x}\| \le \tau\}$, we know from (22) that if $(a,b)$ is sufficiently close to $(0,0)$, there exists $\hat{x} \in \Phi_\tau(a,b)$ such that $\|x - \hat{x}\| \le \kappa\|(a,b)\|$. Therefore, we have

$$\|\hat{x} - \bar{x}\| \le \|\hat{x} - x\| + \|x - \bar{x}\| \le \kappa\|(a,b)\| + \tau,$$

which implies $\hat{x} \in \Phi_\tau(a,b)$. By noting that $\hat{x} \to x$ as $(a,b) \to (0,0)$, we conclude that the set-valued mapping $\Phi_\tau$ is also lower semi-continuous at $(0,0)$.

Next, for any $(a,b)$ sufficiently close to $(0,0)$, consider the following optimization problem

$$\begin{aligned} \min \quad & f(x) - \langle a, x \rangle \\ \text{s.t.} \quad & x \in \Phi_\tau(a,b). \end{aligned} \tag{23}$$

It follows from the Berge theorem (cf. [12, Chapter 9, Theorem 3]) that the solution mapping $X_\tau$ of (23) is nonempty around $(0,0)$ and upper semi-continuous at $(0,0)$. Thus, by noting that $X_\tau(0,0) = X(0,0) = \{\bar{x}\}$, we have for any given $0 < \varepsilon_2' < \tau$ there exists $0 < \varepsilon_1' < \varepsilon_1$ such that for any $\|(a,b)\| < \varepsilon_1'$, $X_\tau(a,b) \subseteq \{x \in \mathcal{X} \mid \|x - \bar{x}\| < \varepsilon_2'\}$, and for any solution $x \in X_\tau(a,b)$ of (23),

$$\|x - \bar{x}\| \le \varepsilon_2' < \tau + \kappa\|(a,b)\|.$$

Therefore, for any $\|(a,b)\| < \varepsilon_1'$, the constraint $\|x - \bar{x}\| \le \tau + \kappa\|(a,b)\|$ of (23) is inactive at any solution point $x \in X_\tau(a,b)$, which implies that

$$\emptyset \neq X_\tau(a,b) \subset X(a,b) \cap \{x \in \mathcal{X} \mid \|x - \bar{x}\| < \varepsilon_2'\} \subset X(a,b) \cap \mathcal{V}.$$

Thus, we know that $X$ is lower semi-continuous at $(0,0,\bar{x}) \in \operatorname{gph} X$. The proof is completed. $\square$

For any given open set $\mathcal{O} \subset \mathcal{X}$, we use $X_{\mathrm{KKT}} \cap \mathcal{O}$ to denote the set-valued mapping $X_{\mathrm{KKT}} \cap \mathcal{O} : \mathcal{X} \times \mathcal{Y} \rightrightarrows \mathcal{X}$ defined by $(X_{\mathrm{KKT}} \cap \mathcal{O})(a,b) = X_{\mathrm{KKT}}(a,b) \cap \mathcal{O}$ for all $(a,b) \in \mathcal{X} \times \mathcal{Y}$. By using the reduction approach, we can extend Robinson's classical result [32, Theorem 2.3] to problem (1) easily. Note that we are not able to apply [32, Theorem 2.3] to problem (1) directly, since the closed convex set $\mathcal{K}$ in problem (1) may not be a cone.

LEMMA 12. *Suppose that $\bar{x}$ is a feasible solution to problem (1) with $(a,b) = (0,0)$. If the RCQ holds at $\bar{x}$, then there exist an open neighborhood $\mathcal{U}$ of the origin and an open neighborhood of $\mathcal{V}$ of $\bar{x}$ such that the Lagrange multipliers mapping $M$ (defined by (5)) is upper semi-continuous on $\mathcal{V} \times \mathcal{U}$ and $X_{\mathrm{KKT}} \cap \mathcal{V}$ is upper semi-continuous on $\mathcal{U}$.*

*Proof.* By the $C^2$-cone reducibility assumption, we know that there exist an open neighborhood $\mathcal{U}$ of the origin and an open neighborhood of $\mathcal{V}$ of $\bar{x}$ such that problem (1) is locally equivalent to its reduced problem (17). By shrinking $\mathcal{U}$ and $\mathcal{V}$ if necessary, we know from (15) of Lemma 8 that for each $(x,a,b) \in \mathcal{V} \times \mathcal{U}$,

$$\mathcal{N}_\mathcal{K}(G(x) + b) = \Xi'(G(x)+b)^* \mathcal{N}_\mathcal{Q}(\mathcal{G}(x,b)).$$

Thus, for each given $(a,b) \in \mathcal{U}$, the sets of stationary points of problems (1) and (17) restricted to $\mathcal{V}$ coincide. Also, we have for each $(x,a,b) \in \mathcal{V} \times \mathcal{U}$,

(24) $$M(x,a,b) = \Xi'(G(x)+b)^* \mathcal{M}(x,a,b),$$

where $\mathcal{M}(x,a,b)$ is the set of Lagrange multipliers of problem (17) associated with $(x,a,b)$. Since for $(a,b) = (0,0)$ the RCQ for problem (1) holds at the feasible point $\bar{x}$ if and only if the RCQ for problem (17) holds at $\bar{x}$, we know from [32, Theorem 2.3] that $\mathcal{M}$ is upper semi-continuous on $\mathcal{V} \times \mathcal{U}$ and $X_{\mathrm{KKT}} \cap \mathcal{V}$ is upper semi-continuous on $\mathcal{U}$. Since $\Xi'$ is continuous on $\mathcal{W}$, we know from (24) that $M$ is also upper semi-continuous on $\mathcal{V} \times \mathcal{U}$. □

By combining Lemmas 11 and 12, we are able to derive the following result on the continuity of the stationary point mapping $X_{\mathrm{KKT}}$, immediately.

COROLLARY 13. *Suppose that $\bar{x}$ is an isolated locally optimal solution of problem (1) with $(a,b) = (0,0)$. If the RCQ (6) holds at $\bar{x}$, then the stationary point mapping $X_{\mathrm{KKT}}$ is lower semi-continuous at $(0,0,\bar{x}) \in \mathrm{gph}\, X_{\mathrm{KKT}}$ and there exists an open neighborhood $\mathcal{V}$ of $\bar{x}$ such that the set-valued mapping $X_{\mathrm{KKT}} \cap \mathcal{V}$ is continuous at $(0,0,\bar{x})$.*

*Proof.* By Lemma 11, we know that for any open neighborhood $\mathcal{V}$ of $\bar{x}$ there exists an open neighborhood $\mathcal{U}$ of the origin such that for any $(a,b) \in \mathcal{U}$, $\emptyset \neq X(a,b) \cap \mathcal{V}$. By shrinking $\mathcal{V}$ and $\mathcal{U}$ if necessary, we know that the RCQ holds at any $x \in X(a,b) \cap \mathcal{V}$ for any $(a,b) \in \mathcal{U}$ (cf. e.g., [7, Remark 2.88]), which implies that $\emptyset \neq X(a,b) \cap \mathcal{V} \subset X_{\mathrm{KKT}}(a,b) \cap \mathcal{V}$. Therefore, we know that $X_{\mathrm{KKT}}$ is lower semi-continuous at $(0,0,\bar{x})$. Finally, it follows from Lemma 12 that there exists an open neighborhood $\mathcal{V}$ such that $X_{\mathrm{KKT}} \cap \mathcal{V}$ is upper semi-continuous at $(0,0)$. Therefore, we conclude that $X_{\mathrm{KKT}} \cap \mathcal{V}$ is continuous at $(0,0,\bar{x})$. □

By Corollary 13 and Lemma 12, we obtain the following result on the lower semi-continuity of $S_{\mathrm{KKT}}$.

PROPOSITION 14. *Suppose that $\bar{x}$ is an isolated locally optimal solution of problem (1) with $(a,b) = (0,0)$ and the corresponding set of Lagrange multipliers $M(\bar{x},0,0) \neq \emptyset$. If the SRCQ (11) holds at $\bar{x}$ with respect to $\bar{y} \in M(\bar{x},0,0)$, then the KKT solution mapping $S_{\mathrm{KKT}}$ is lower semi-continuous at $(0,0,\bar{x},\bar{y}) \in \mathrm{gph}\, S_{\mathrm{KKT}}$.*

*Proof.* Since the SRCQ holds at $\bar{x}$ with respect to $\bar{y}$, we know that the RCQ holds at $\bar{x}$ with respect to $(a,b) = (0,0)$. By Corollary 13, we know that $X_{KKT}$ is lower semi-continuous at $(0,0,\bar{x}) \in \mathrm{gph}\, X_{\mathrm{KKT}}$, i.e., for any open neighborhood $\mathcal{V}_1$ of $\bar{x}$, there exists an open neighborhood $\mathcal{U}$ of the origin such that for every $(a,b) \in \mathcal{U}$, $X_{\mathrm{KKT}}(a,b) \cap \mathcal{V}_1 \neq \emptyset$. Furthermore, by Lemma 12, we know that $M$ is locally bounded at $(\bar{x},0,0)$, i.e., there exists a bounded set $\mathcal{B} \subset \mathcal{Y}$ such that for any $(x,a,b) \in \mathcal{V}_1 \times \mathcal{U}$, $M(x,a,b) \subset \mathcal{B}$.

Next, we shall show that for any open neighborhood $\mathcal{V}_1$ of $\bar{x}$, the open neighborhood $\mathcal{U}$ of the origin (shrinking if necessary) satisfies that for any open neighborhood $\mathcal{V}_2$ of $\bar{y}$ in $\mathcal{Y}$, $M(x,a,b) \cap \mathcal{V}_2 \neq \emptyset$ for every $(x,a,b) \in \mathcal{V}_1 \times \mathcal{U}$. Assume, on the contrary, that there exist an open neighborhood $\widetilde{\mathcal{V}}_2$ of $\bar{y}$ and the sequence $\{(x^k, a^k, b^k)\}$ with $(a^k, b^k) \to (0,0)$ and $X_{\text{KKT}}(a^k, b^k) \ni x^k \to \bar{x}$ such that $M(x^k, a^k, b^k) \cap \widetilde{\mathcal{V}}_2 = \emptyset$, i.e., $y \notin \widetilde{\mathcal{V}}_2$ for any $y \in M(x^k, a^k, b^k)$. Thus, we know that there exists a positive number $\eta > 0$ such that for all $(x^k, a^k, b^k)$,

$$\|y - \bar{y}\| \geq \eta \quad \forall\, y \in M(x^k, a^k, b^k).$$

Since $M$ is locally bounded at $(\bar{x}, 0, 0)$, we may take a bounded sequence $\{y^k\}$ with $y^k \in M(x^k, a^k, b^k)$ for each $k$, which has an accumulation point $\hat{y} \neq \bar{y}$. Since the support function $\sigma(\cdot, \mathcal{K})$ of $\mathcal{K}$ is a closed proper convex function, it follows from [33, Theorem 24.4] that the graph of $\partial\sigma(\cdot, \mathcal{K})$ is a closed subset in $\mathcal{Y} \times \mathcal{Y}$. Therefore, by taking the limit with $k$ in the KKT system (4), we know that $\hat{y} \in M(\bar{x}, 0, 0)$, which contradicts with the fact that $M(\bar{x}, 0, 0)$ is a singleton. Thus, for any open neighborhood $\mathcal{V} = \mathcal{V}_1 \times \mathcal{V}_2$ of $(\bar{x}, \bar{y})$, there exists an open neighborhood $\mathcal{U}$ of the origin such that for any $(a, b) \in \mathcal{U}$, $S_{\text{KKT}}(a, b) \cap \mathcal{V} \neq \emptyset$. This completes the proof. □

By using [32, Theorem 2.4] and the reduction approach, we are able to derive the following proposition, in which a sufficient condition for a locally optimal solution $\bar{x}$ to be isolated is established.

PROPOSITION 15. *Suppose that the RCQ holds at a locally optimal solution $\bar{x}$ of problem (1) with $(a, b) = (0, 0)$ and that Robinson's SOSC*

(25) $$\inf_{y \in M(\bar{x}, 0, 0)} \left\{ \langle d, \nabla^2_{xx} L(\bar{x}; y) d \rangle - \sigma\left( y, \mathcal{T}^2_{\mathcal{K}}(G(\bar{x}), G'(\bar{x})d) \right) \right\} > 0 \quad \forall\, d \in \mathcal{C}(\bar{x}) \setminus \{0\}$$

*holds at $\bar{x}$. Then, there exists an open neighborhood $\mathcal{V}$ of $\bar{x}$ such that $X_{\text{KKT}}(0, 0) \cap \mathcal{V} = \{x\}$, which implies that $\bar{x}$ is an isolated locally optimal solution of problem (1) with $(a, b) = (0, 0)$.*

*Proof.* By the $C^2$-cone reducibility assumption, we know that there exists an open neighborhood $\mathcal{V}$ of $\bar{x}$ such that problem (1) with $(a, b) = (0, 0)$ is locally equivalent to its reduced problem (17). For $(a, b) = (0, 0)$, the RCQ for problem (1) holds at $\bar{x}$ if and only if the RCQ for problem (17) holds at $\bar{x}$. Let $\mathcal{M}(\bar{x}, 0, 0)$ be the set of Lagrange multipliers of problem (17) associated with $\bar{x}$ and $(a, b) = (0, 0)$. It follows from Lemma 8 that for any $y \in M(\bar{x}, 0, 0)$, there is a unique $u \in \mathcal{M}(\bar{x}, 0, 0)$ such that $y = \Xi'(G(\bar{x}))^* u$, i.e., $\Xi'(G(\bar{x}))^*$ is a bijection between $\mathcal{M}(\bar{x}, 0, 0)$ and $M(\bar{x}, 0, 0)$. By [7, (3.271)], we know that the corresponding critical cone of problem (17) with $(a, b) = (0, 0)$ coincides with the critical cone $\mathcal{C}(\bar{x})$ of problem (1) defined by (9). Thus, we know from [32, Theorem 2.4] that if

(26) $$\inf_{u \in \mathcal{M}(\bar{x}, 0, 0)} \langle d, \nabla^2_{xx} \mathcal{L}(\bar{x}; u) d \rangle > 0 \quad \forall\, d \in \mathcal{C}(\bar{x}) \setminus \{0\},$$

then $\bar{x}$ is the unique stationary point of problem (17) with $(a, b) = (0, 0)$ in $\mathcal{V}$, where $\mathcal{L}(x; u) : \mathcal{V} \times \mathcal{Z} \to \Re$ is the Lagrangian function of problem (17) defined by

$$\mathcal{L}(x; u) := f(x) + \langle u, \Xi(G(x)) \rangle, \quad (x, u) \in \mathcal{V} \times \mathcal{Z}.$$

By shrinking $\mathcal{V}$ if necessary, we know from (15) of Lemma 8 that the sets of stationary points of problems (1) and (17) with $(a, b) = (0, 0)$ restricted to $\mathcal{V}$ are the same. Thus, we know that $X_{\mathrm{KKT}}(0,0) \cap \mathcal{V} = \{\bar{x}\}$. Note that by the chain rule, we have for any $u \in \mathcal{M}(\bar{x}, 0, 0)$ and $d \in \mathcal{C}(\bar{x})$,

$$\langle d, \nabla^2_{xx}\mathcal{L}(\bar{x}; u)d \rangle = \langle d, \nabla^2_{xx}L(\bar{x}; y)d \rangle + \langle u, \Xi''(G(\bar{x}))(G'(\bar{x})d, G'(\bar{x})d) \rangle.$$

Thus, it follows from (15) of Lemma 8 that the conditions (25) and (26) are equivalent. This completes the proof. □

REMARK 2. *The quadratic growth condition at $\bar{x}$ for problem (1) with $(a, b) = (0, 0)$, which implies that $\bar{x}$ is a strict locally optimal solution, may not guarantee that $\bar{x}$ is an isolated local minimizer. For instance, consider the following example [32]*

$$\begin{aligned} \min \quad & x^2/2 \\ \text{s.t.} \quad & x^6 \sin(1/x) = 0, \end{aligned}$$

*where $x \in \Re$ and $\sin(1/0) := 0$. It is clear that the quadratic growth condition holds at the strict locally optimal solution $\bar{x} = 0$. However, each element of the set $\{(n\pi)^{-1} \mid n = \pm 1, \pm 2, \ldots\}$ is a locally optimal solution and $\bar{x} = 0$ is its cluster point, which implies that $\bar{x} = 0$ is not an isolated locally optimal solution.*

By combining Proposition 15 and Corollary 13, we obtain the following result on the continuity of $X_{\mathrm{KKT}}$, immediately.

THEOREM 16. *Suppose that the RCQ holds at a locally optimal solution $\bar{x}$ of problem (1) with $(a, b) = (0, 0)$ and Robinson's SOSC (25) holds at $\bar{x}$. Then there exists an open neighborhood $\mathcal{V}$ of $\bar{x}$ such that the set-valued mapping $X_{\mathrm{KKT}} \cap \mathcal{V}$ is continuous at $(0, 0, \bar{x})$.*

Note that if the set of Lagrange multipliers with respect to a locally optimal solution is a singleton, then Robinson's SOSC (25) coincides with the SOSC defined by (14). Next, we shall provide an example to illustrate that the uniqueness of Lagrange multipliers is not a necessary condition for the continuity of $X_{\mathrm{KKT}}$.

EXAMPLE 2. *Consider the following convex quadratic SDP problem:*

(27)
$$\begin{aligned} (P) \quad \min \quad & \tfrac{1}{2}x^2 - x + t \\ \text{s.t.} \quad & tA - xI \in \mathcal{S}^2_+, \\ & t \geq 0 \end{aligned}$$

*and its dual*

(28)
$$\begin{aligned} (D) \quad \max \quad & -\tfrac{1}{2}(\langle I, Y \rangle - 1)^2 \\ \text{s.t.} \quad & \langle A, Y \rangle \leq 1, \\ & Y \in \mathcal{S}^2_+, \end{aligned}$$

*where $x, t \in \Re$, $A = \begin{bmatrix} 1 & -2 \\ -2 & 1 \end{bmatrix}$ and $I$ is the $2 \times 2$ identity matrix. We know that $(\bar{x}, \bar{t}) = (0, 0)$ is the unique optimal solution to the primal problem (27), and the RCQ holds at $(\bar{x}, \bar{t})$. Since the critical cone of problem (27) at $(\bar{x}, \bar{t})$ is $\mathcal{C}(\bar{x}, \bar{t}) = \{0\}$, we know that Robinson's SOSC holds at $(\bar{x}, \bar{t})$. Thus, by Theorem 16, we know that $X_{\mathrm{KKT}} \cap \mathcal{V}$ is continuous at $(0, 0, \bar{x}, \bar{t})$. Note that the solution set to the dual problem (28) is given by*

$$\{Y \in \mathcal{S}^2_+ \mid \langle A, Y \rangle \leq 1, \ \langle I, Y \rangle = 1\}.$$

*This implies that the Lagrange multipliers of problem (27) are not unique.*

Since the SRCQ implies the uniqueness of Lagrange multipliers, by combining Propositions 14 and 15, we obtain the following result, directly.

THEOREM 17. *Let $\bar{x}$ be a feasible solution of problem (1) with $(a,b) = (0,0)$. Suppose that the SRCQ holds at $\bar{x}$ with respect to $\bar{y} \in M(\bar{x}, 0, 0) \neq \emptyset$ and the SOSC (14) holds at $\bar{x}$ for problem (1) with respect to $(a,b) = (0,0)$. Then the set-valued mapping $S_{\mathrm{KKT}}$ is lower semi-continuous at $(0, 0, \bar{x}, \bar{y}) \in \mathrm{gph}\, S_{\mathrm{KKT}}$.*

**4. The robust isolated calmness of $S_{\mathrm{KKT}}$.** For each given $(a,b) \in \mathcal{X} \times \mathcal{Y}$, it is well-known (e.g., [13, Proposition 1.5.9]) that the set of solutions of the KKT system (4) can be rewritten as

$$\text{(29)} \qquad S_{\mathrm{KKT}}(a,b) = \big\{(x, z - \Pi_{\mathcal{K}}(z)) \in \mathcal{X} \times \mathcal{Y} \mid \Psi(x, z) = (a, -b)\big\},$$

where $\Psi : \mathcal{X} \times \mathcal{Y} \to \mathcal{X} \times \mathcal{Y}$ is Robinson's normal mapping defined by

$$\text{(30)} \qquad \Psi(x, z) = \begin{bmatrix} \nabla f(x) + G'(x)^*(z - \Pi_{\mathcal{K}}(z)) \\ G(x) - \Pi_{\mathcal{K}}(z) \end{bmatrix}, \quad (x, z) \in \mathcal{X} \times \mathcal{Y}.$$

Let $(\bar{x}, \bar{y})$ be a solution to the KKT system (4) with $(a,b) = (0,0)$. Denote $\bar{z} := G(\bar{x}) + \bar{y}$. Then, since $\Pi_{\mathcal{K}}$ is globally Lipschitz continuous (with modulus 1) and $G$ is locally Lipschitz continuous, it is an easy exercise to see that the KKT solution mapping $S_{\mathrm{KKT}}$ is isolated calm (robustly isolated calm, respectively) at the origin for $(\bar{x}, \bar{y})$ if and only if the set-valued mapping $\Psi^{-1}$ is isolated calm (robustly isolated calm, respectively) at the origin for $(\bar{x}, \bar{z})$. Moreover, $S_{\mathrm{KKT}}$ has the Aubin property at the origin for $(\bar{x}, \bar{y})$ if and only if $\Psi^{-1}$ has the Aubin property at the origin for $(\bar{x}, \bar{z})$.

When $(a,b) = (0,0)$, the KKT system (4) is equivalent to the following system of nonsmooth equations:

$$\text{(31)} \qquad\qquad\qquad F(x,y) = 0,$$

where $F : \mathcal{X} \times \mathcal{Y} \to \mathcal{X} \times \mathcal{Y}$ is the natural mapping defined by

$$\text{(32)} \qquad F(x,y) := \begin{bmatrix} \nabla f(x) + G'(x)^* y \\ G(x) - \Pi_{\mathcal{K}}(G(x) + y) \end{bmatrix}, \quad (x, y) \in \mathcal{X} \times \mathcal{Y}.$$

It is clear that $(0, 0, \bar{x}, \bar{y}) \in \mathrm{gph}\, S_{\mathrm{KKT}}$ if and only if $(0, 0, \bar{x}, \bar{y}) \in \mathrm{gph}\, F^{-1}$. The following simple observation can be derived from the definition of the isolated calmness without assuming that $\mathcal{K}$ is $C^2$-cone reducible.

LEMMA 18. *Let $(0, 0, \bar{x}, \bar{y}) \in \mathrm{gph}\, S_{\mathrm{KKT}}$. The set-valued mapping $S_{\mathrm{KKT}}$ is isolated calm at the origin for $(\bar{x}, \bar{y})$ if and only if the set-valued mapping $F^{-1}$ is isolated calm at the origin for $(\bar{x}, \bar{y})$.*

*Proof.* Firstly, suppose that $S_{\mathrm{KKT}}$ is isolated calm at the origin for $(\bar{x}, \bar{y})$. Then, we know that there exist $\kappa > 0$, $\varepsilon_1 > 0$ and $\varepsilon_2 > 0$ such that for any $(a,b) \in \mathcal{X} \times \mathcal{Y}$ satisfying $\|(a,b)\| < \varepsilon_1$,

$$\text{(33)} \quad S_{\mathrm{KKT}}(a,b) \cap \{(x,y) \in \mathcal{X} \times \mathcal{Y} \mid \|(x,y) - (\bar{x},\bar{y})\| < \varepsilon_2\} \subset \{(\bar{x},\bar{y})\} + \kappa \|(a,b)\| \mathbb{B}_{\mathcal{X} \times \mathcal{Y}}.$$

By noting that $G$ is continuously differentiable, we know that there exists a constant $\varepsilon_3 > 0$ such that $\|G'(x)\| \leq \varepsilon_3$ for any $(x,y)$ satisfying $\|(x,y) - (\bar{x},\bar{y})\| < \varepsilon_2$. Choose $0 < \eta_1 < \min\{\varepsilon_1/\sqrt{2\varepsilon_3^2 + 2}, \varepsilon_2/2\}$ and $0 < \eta_2 < \varepsilon_2/2$. For any $(\hat{a}, \hat{b}) \in \mathcal{X} \times \mathcal{Y}$ satisfying $\|(\hat{a}, \hat{b})\| < \eta_1$, let $(\hat{x}, \hat{y}) \in F^{-1}(\hat{a}, \hat{b}) \cap \{(x,y) \mid \|(x,y) - (\bar{x},\bar{y})\| < \eta_2\}$ be arbitrarily given. Therefore, it follows from (31) that

$$\begin{cases} \nabla f(\hat{x}) + G'(\hat{x})^* \hat{y} = \hat{a}, \\ G(\hat{x}) - \hat{b} - \Pi_{\mathcal{K}}(G(\hat{x}) - \hat{b} + \hat{b} + \hat{y}) = 0, \end{cases}$$

which is equivalent to
$$\begin{cases} \nabla f(\hat{x}) + G'(\hat{x})^*(\hat{y} + \hat{b}) = \hat{a} + G'(\hat{x})^*\hat{b}, \\ \hat{y} + \hat{b} \in \mathcal{N}_{\mathcal{K}}(G(\hat{x}) - \hat{b}). \end{cases}$$

Thus, by (4), we have $(\hat{x}, \hat{y} + \hat{b}) \in S_{\text{KKT}}(\hat{a} + G'(\hat{x})^*\hat{b}, -\hat{b})$. Moreover, since $\|(\hat{x}, \hat{y}) - (\bar{x}, \bar{y})\| < \eta_2$, we have
$$\|(\hat{x}, \hat{y} + \hat{b}) - (\bar{x}, \bar{y})\| \leq \|(\hat{x}, \hat{y}) - (\bar{x}, \bar{y})\| + \|\hat{b}\| < \eta_2 + \|(\hat{a}, \hat{b})\| < \frac{\varepsilon_2}{2} + \eta_1 < \varepsilon_2,$$

which implies that $(\hat{x}, \hat{y} + \hat{b}) \in S_{\text{KKT}}(\hat{a} + G'(\hat{x})^*\hat{b}, -\hat{b}) \cap \{(x, y) \mid \|(x, y) - (\bar{x}, \bar{y})\| < \varepsilon_2\}$. By noting that $\|(\hat{a} + G'(\hat{x})^*\hat{b}, -\hat{b})\| \leq \sqrt{2\varepsilon_3^2 + 2}\|(\hat{a}, \hat{b})\| < \varepsilon_1$, we know from (33) that
$$\|(\hat{x}, \hat{y} + \hat{b}) - (\bar{x}, \bar{y})\| \leq \kappa\sqrt{2\varepsilon_3^2 + 2}\|(\hat{a}, \hat{b})\|.$$

Thus, we have
$$\|\hat{x} - \bar{x}\|^2 + \frac{1}{2}\|\hat{y} - \bar{y}\|^2 - \|\hat{b}\|^2 \leq \|\hat{x} - \bar{x}\|^2 + \|\hat{y} - \bar{y} + \hat{b}\|^2 \leq \kappa^2(2\varepsilon_3^2 + 2)\|(\hat{a}, \hat{b})\|^2.$$

Hence, $\|\hat{x} - \bar{x}\|^2 + \|\hat{y} - \bar{y}\|^2 \leq 4(\kappa^2(\varepsilon_3^2 + 1) + 1)\|(\hat{a}, \hat{b})\|^2$. Since $(\hat{x}, \hat{y})$ is arbitrary, we obtain that for any $\|(\hat{a}, \hat{b})\| < \eta_1$,

(34) $\quad F^{-1}(\hat{a}, \hat{b}) \cap \{(x, y) \mid \|(x, y) - (\bar{x}, \bar{y})\| < \eta_2\} \subset \{(\bar{x}, \bar{y})\} + \tau\|(\hat{a}, \hat{b})\|\mathbb{B}_{\mathcal{X} \times \mathcal{Y}},$

where $\tau = 2\sqrt{\kappa^2(\varepsilon_3^2 + 1) + 1} > 0$. Thus, $F^{-1}$ is isolated calm at the origin for $(\bar{x}, \bar{y})$.

Conversely, assume that $F^{-1}$ is isolated calm at the origin for $(\bar{x}, \bar{y})$. Then, there exist $\tau > 0$, $\eta_1 > 0$ and $\eta_2 > 0$ such that for any $\|(\hat{a}, \hat{b})\| < \eta_1$, (34) holds. Again, by the continuity of $G'(x)$, we know that there exists a constant $\eta_3 > 0$ such that for any $\|(x, y) - (\bar{x}, \bar{y})\| < \eta_2$, $\|G'(x)\| \leq \eta_3$. Choose $0 < \varepsilon_1 < \min\{\eta_1/\sqrt{2\eta_3^2 + 2}, \eta_2/2\}$ and $0 < \varepsilon_2 < \eta_2/2$. For any $(a, b) \in \mathcal{X} \times \mathcal{Y}$ satisfying $\|(a, b)\| < \varepsilon_1$, let $(\tilde{x}, \tilde{y}) \in S_{\text{KKT}}(a, b) \cap \{(x, y) \mid \|(x, y) - (\bar{x}, \bar{y})\| < \varepsilon_2\}$ be arbitrarily given. We know from (4) that
$$\begin{cases} \nabla f(\tilde{x}) + G'(\tilde{x})^*(\tilde{y} + b) = a + G'(\tilde{x})^*b, \\ G(\tilde{x}) - \Pi_{\mathcal{K}}(G(\tilde{x}) + b + \tilde{y}) = -b. \end{cases}$$

Thus, we know from (31) that $(\tilde{x}, \tilde{y} + b) \in F^{-1}(a + G'(\tilde{x})^*b, -b)$. Moreover, since $\|(\tilde{x}, \tilde{y}) - (\bar{x}, \bar{y})\| < \varepsilon_2$, we obtain that
$$\|(\tilde{x}, \tilde{y} + b) - (\bar{x}, \bar{y})\| \leq \|(\tilde{x}, \tilde{y}) - (\bar{x}, \bar{y})\| + \|b\| < \varepsilon_2 + \|(a, b)\| < \frac{\eta_2}{2} + \varepsilon_1 < \eta_2,$$

which implies that $(\tilde{x}, \tilde{y} + b) \in F^{-1}(a + G'(\tilde{x})^*b, -b) \cap \{(x, y) \mid \|(x, y) - (\bar{x}, \bar{y})\| < \eta_2\}$. Again, by noting that $\|(a + G'(\tilde{x})^*b, -b)\| \leq \sqrt{2\eta_3^2 + 2}\|(a, b)\| < \eta_1$, we know from (34) that
$$\|\tilde{x} - \bar{x}\|^2 + \frac{1}{2}\|\tilde{y} - \bar{y}\|^2 - \|b\|^2 \leq \|\tilde{x} - \bar{x}\|^2 + \|\tilde{y} - \bar{y} + b\|^2 \leq \tau^2(2\eta_3^2 + 2)\|(a, b)\|^2.$$

Therefore, $\|\tilde{x} - \bar{x}\|^2 + \|\tilde{y} - \bar{y}\|^2 \leq 4(\tau^2(\eta_3^2 + 1) + 1)\|(a, b)\|^2$. Since $(\tilde{x}, \tilde{y})$ is arbitrary, we obtain that for any $\|(a, b)\| < \varepsilon_1$, (33) holds with $\kappa = 2\sqrt{\tau^2(\eta_3^2 + 1) + 1} > 0$. Thus, $S_{\text{KKT}}$ is isolated calm at the origin for $(\bar{x}, \bar{y})$. This completes the proof. □

Under the assumption that $\mathcal{K}$ is $C^2$-cone reducible, $F$ is locally Lipschitz continuous around $(\bar{x}, \bar{y})$ and is directionally differentiable at $(\bar{x}, \bar{y}) \in \mathcal{X} \times \mathcal{Y}$, the following result of the basic characterization on the isolated calmness of $F^{-1}$ follows from Lemma 4 and [34, 8(19)], directly.

LEMMA 19. *Let $\bar{x}$ be a stationary point of problem (1) with $(a, b) = (0, 0)$. Suppose that $\bar{y} \in M(\bar{x}, 0, 0)$. Then, the set-valued mapping $F^{-1}$ is isolated calm at the origin for $(\bar{x}, \bar{y})$ if and only if $F'((\bar{x}, \bar{y}); (\triangle x, \triangle y)) = 0$ implies $(\triangle x, \triangle y) = 0$, i.e.,*

$$
(35) \quad \begin{cases} \nabla_{xx}^2 L(\bar{x}; \bar{y})\triangle x + G'(\bar{x})^* \triangle y = 0, \\ G'(\bar{x})\triangle x - \Pi'_{\mathcal{K}}(G(\bar{x}) + \bar{y}; G'(\bar{x})\triangle x + \triangle y) = 0 \end{cases} \implies (\triangle x, \triangle y) = 0.
$$

REMARK 3. *Note that both the normal mapping (30) and the natural mapping (32) are locally Lipschitz continuous and directionally differentiable under the assumption that $\mathcal{K}$ is $C^2$-cone reducible. It is easy to check that (35) is equivalent to $\Psi'((\bar{x}, \bar{z}); (\triangle x, \triangle z)) = 0$ implies $(\triangle x, \triangle z) = 0$, i.e., $\Psi^{-1}$ is isolated calmness at the origin for $(\bar{x}, \bar{z})$, where $\bar{z} = G(\bar{x}) + \bar{y}$. Therefore, under the $C^2$-cone reducibility assumption, one may derive Lemma 18 by employing Lemma 19 directly.*

REMARK 4. *Recently, under the assumption that the constraint non-degeneracy (see (12) for the definition) holds at a stationary point, Mordukhovich et al. [26, 27] established a characterization on the isolated calmness of the solution map of the parameter-dependent generalized equation. Moreover, when applied to the canonically perturbed problem (1), the characterization proposed in [26, 27] becomes (35). However, the constraint non-degeneracy assumption seems to be crucial for the analysis conducted in [26, 27], in particular [26, Lemma 3.1] and [27, Proposition 3.1].*

It is known from Klatte and Kummer [20], generalizing a result of Fusek [15] on the nonlinear SDP problem, that for problem (1) ($\mathcal{K}$ does not need to be $C^2$-cone reducible) with $(a, b) = (0, 0)$, if $S_{\mathrm{KKT}}$ has the Aubin property at the origin for $(\bar{x}, \bar{y})$, then the constraint non-degeneracy condition (12) holds at $\bar{x}$. Therefore, by using the fact that for a $C^2$-cone reducible convex set $\mathcal{K}$, $\Pi_{\mathcal{K}}$ is directionally differentiable at any point $y \in \mathcal{Y}$ [3, Theorem 7.2], we obtain from Lemma 5 the following result on the relationship between the Aubin property[2] of $F^{-1}$ and the isolated calmness of $F^{-1}$.

PROPOSITION 20. *Let $\bar{x}$ be a stationary point for problem (1) with $(a, b) = (0, 0)$. Suppose that $F^{-1}$ has the Aubin property at the origin for $(\bar{x}, \bar{y})$ with $y \in M(\bar{x}, 0, 0) \neq \emptyset$, then the constraint non-degeneracy condition (12) holds at $\bar{x}$ and $F^{-1}$ is isolated calm at the origin for $(\bar{x}, \bar{y})$.*

Next, under the assumption that $\mathcal{K}$ is $C^2$-cone reducible, by employing Lemma 10 we are able to extend the result in [16, Theorem 4.1 (i)] on the isolated calmness of $F^{-1}$.

PROPOSITION 21. *Let $\bar{x}$ be a stationary point for problem (1) with $(a, b) = (0, 0)$. If the SOSC (14) for problem (1) holds at $\bar{x}$ with respect to $(a, b) = (0, 0)$ and the SRCQ (11) holds at $\bar{x}$ with respect to $\bar{y} \in M(\bar{x}, 0, 0)$, then $F^{-1}$ is isolated calm at the origin for $(\bar{x}, \bar{y})$.*

*Proof.* Let $(\triangle x, \triangle y) \in \mathcal{X} \times \mathcal{Y}$ be arbitrarily chosen such that

$$
(36) \quad \begin{cases} \nabla_{xx}^2 L(\bar{x}; \bar{y})\triangle x + G'(\bar{x})^* \triangle y = 0, \\ G'(\bar{x})\triangle x - \Pi'_{\mathcal{K}}(G(\bar{x}) + \bar{y}; G'(\bar{x})\triangle x + \triangle y) = 0. \end{cases}
$$

---

[2]The above mentioned result of Klatte and Kummer [20] on the constraint non-degeneracy condition can be easily applied to $F^{-1}$ if $\mathcal{K}$ is a $C^2$-cone reducible convex set.

By part (i) of Lemma 10, we know from the second equation of (36) that

$$G'(\bar{x})\triangle x \in \mathcal{C}_{\mathcal{K}}(G(\bar{x}), \bar{y}) \quad \text{and} \quad \langle G'(\bar{x})\triangle x, \triangle y\rangle = -\sigma(\bar{y}, \mathcal{T}_{\mathcal{K}}^2(G(\bar{x}), G'(\bar{x})\triangle x)).$$

Thus, we have $\triangle x \in \mathcal{C}(\bar{x})$. By taking the inner product between $\triangle x$ and both sides of the first equation of (36), respectively, we obtain that

$$\begin{aligned} 0 &= \langle \triangle x, \nabla_{xx}^2 L(\bar{x}; \bar{y})\triangle x\rangle + \langle G'(\bar{x})\triangle x, \triangle y\rangle \\ &= \langle \triangle x, \nabla_{xx}^2 L(\bar{x}; \bar{y})\triangle x\rangle - \sigma(\bar{y}, \mathcal{T}_{\mathcal{K}}^2(G(\bar{x}), G'(\bar{x})\triangle x)). \end{aligned}$$

Hence, it follows from the SOSC (14) for problem (1) with respect to $(a, b) = (0, 0)$ that $\triangle x = 0$. Therefore, (36) is reduced to

$$\begin{cases} G'(\bar{x})^* \triangle y = 0, \\ \Pi'_{\mathcal{K}}(G(\bar{x}) + \bar{y}; \triangle y) = 0. \end{cases}$$

By part (ii) of Lemma 10, we have

$$\triangle y \in \left[G'(\bar{x})\mathcal{X} + \mathcal{T}_{\mathcal{K}}(G(\bar{x})) \cap \bar{y}^{\perp}\right]^{\circ}.$$

Thus, we know from the SRCQ (11) that $\triangle y = 0$. Therefore, the only $(\triangle x, \triangle y) \in \mathcal{X} \times \mathcal{Y}$ satisfying (36) is $(\triangle x, \triangle y) = 0$. By Lemma 19, we know that $F^{-1}$ is isolated calm at the origin for $(\bar{x}, \bar{y})$. This completes the proof. □

In the following, we shall study the converse implication of Proposition 21. Firstly, we introduce the following result.

LEMMA 22. *Suppose that $\bar{x}$ is a stationary point to problem (1) with $(a, b) = (0, 0)$ and the SRCQ (11) holds at $\bar{x}$ with respect to $\bar{y} \in M(\bar{x}, 0, 0)$. Suppose that $F^{-1}$ is isolated calm at the origin for $(\bar{x}, \bar{y})$ and that there exists $\triangle x \in \mathcal{C}(\bar{x}) \setminus \{0\}$ such that*

$$(37) \qquad \langle \triangle x, \nabla_{xx}^2 L(\bar{x}; \bar{y})\triangle x\rangle + \Upsilon(G'(\bar{x})\triangle x) = 0,$$

*where $\Upsilon(\cdot)$ is the quadratic function defined by (19). Then there exists $\bar{d} \in \mathcal{C}(\bar{x})$ such that*

$$(38) \qquad \langle \nabla_{xx}^2 L(\bar{x}; \bar{y})\triangle x, \bar{d}\rangle + \frac{1}{2}\langle G'(\bar{x})^* \nabla \Upsilon(G'(\bar{x})\triangle x), \bar{d}\rangle < 0,$$

*Proof.* Suppose that for every $d \in \mathcal{C}(\bar{x})$, the inequality (38) fails to hold. By (37) and noting that $\langle \nabla \Upsilon(G'(\bar{x})\triangle x), G'(\bar{x})\triangle x\rangle = 2\Upsilon(G'(\bar{x})\triangle x)$, we know that $\triangle x$ is an optimal solution to the following linear conic programming problem

$$(39) \qquad \begin{array}{ll} \min & \langle \nabla_{xx}^2 L(\bar{x}; \bar{y})\triangle x, d\rangle + \frac{1}{2}\langle G'(\bar{x})^* \nabla \Upsilon(G'(\bar{x})\triangle x), d\rangle \\ \text{s.t.} & G'(\bar{x})d \in \mathcal{C}_{\mathcal{K}}(G(\bar{x}), \bar{y}). \end{array}$$

It is clear that the RCQ of problem (39) holds at $\triangle x$, since the SRCQ (11) holds at $\bar{x}$ with respect to $\bar{y}$. Therefore, we know that there exists $\triangle \eta \in \mathcal{Y}$ such that

$$(40) \qquad \begin{cases} \nabla_{xx}^2 L(\bar{x}; \bar{y})\triangle x + \frac{1}{2}G'(\bar{x})^* \nabla \Upsilon(G'(\bar{x})\triangle x) + G'(\bar{x})^* \triangle \eta = 0, \\ G'(\bar{x})\triangle x \in \mathcal{C}_{\mathcal{K}}(G(\bar{x}), \bar{y}), \\ \triangle \eta \in \mathcal{N}_{\mathcal{C}_{\mathcal{K}}(G(\bar{x}), \bar{y})}(G'(\bar{x})\triangle x). \end{cases}$$

Denote $\triangle y := \triangle \eta + \frac{1}{2}\nabla \Upsilon(G'(\bar{x})\triangle x)$. By noting that $\triangle x \in \mathcal{C}(\bar{x})$ and

$$\langle G'(\bar{x})\triangle x, \nabla \Upsilon(G'(\bar{x})\triangle x)\rangle = 2\Upsilon(G'(\bar{x})\triangle x),$$

we obtain that
$$\langle G'(\bar{x})\triangle x, \triangle y\rangle = \Upsilon(G'(\bar{x})\triangle x) = -\sigma(\bar{y}, \mathcal{T}_\mathcal{K}^2(G(\bar{x}), G'(\bar{x})\triangle x)).$$

Therefore, it follows from part (i) of Lemma 10 that
$$G'(\bar{x})\triangle x - \Pi'_\mathcal{K}(G(\bar{x}) + \bar{y}; G'(\bar{x})\triangle x + \triangle y) = 0,$$

which, together with the first equations of (40), implies that $0 \neq (\triangle x, \triangle y)$ satisfies the equations in (35). This contradicts the isolated calmness of $F^{-1}$ at origin for $(\bar{x}, \bar{y})$. This contradiction completes the proof. $\square$

By using Lemma 22, we show in the following proposition that the converse implication of Proposition 21 also holds.

PROPOSITION 23. *Suppose that $\bar{x}$ is a locally optimal solution to problem* (1) *with $(a,b) = (0,0)$ and the RCQ holds at $\bar{x}$. If $F^{-1}$ is isolated calm at the origin for $(\bar{x}, \bar{y})$, then*

(i) *the SRCQ* (11) *holds at $\bar{x}$ with respect to $\bar{y}$;*

(ii) *the SOSC* (14) *for problem* (1) *with $(a,b) = (0,0)$ holds at $\bar{x}$.*

*Proof.* We first prove part (i). Assume, on the contrary, that the SRCQ (11) does not hold at $\bar{x}$ with respect to $\bar{y}$. Then, there exists $0 \neq \triangle y \in \mathcal{Y}$ such that
$$\triangle y \in \left[G'(\bar{x})\mathcal{X} + \mathcal{T}_\mathcal{K}(G(\bar{x})) \cap \bar{y}^\perp\right]^\circ.$$

Thus, we know from part (ii) of Lemma 10 that
$$\begin{cases} G'(\bar{x})^* \triangle y = 0, \\ \Pi'_\mathcal{K}(G(\bar{x}) + \bar{y}; \triangle y) = 0, \end{cases}$$

which implies that $F'((\bar{x}, \bar{y}); (0, \triangle y)) = 0$. Since $F^{-1}$ is isolated calm at the origin for $(\bar{x}, \bar{y})$, we know from Lemma 19 that $\triangle y = 0$. This contradiction shows that the SRCQ (11) holds at $\bar{x}$ with respect to $\bar{y}$.

Next, we shall prove part (ii). Since $\bar{x}$ is a locally optimal solution to problem (1) with $(a,b) = (0,0)$ and the SRCQ holds at $\bar{x}$, we know from Theorem 7 that

$$\langle d, \nabla^2_{xx} L(\bar{x}; \bar{y}) d \rangle + \Upsilon(G'(\bar{x})d) \geq 0 \quad \forall d \in \mathcal{C}(\bar{x}), \tag{41}$$

where $\Upsilon(\cdot)$ is the quadratic function defined by (19). Therefore, if we assume that the SOSC (14) does not hold at $\bar{x}$, then there exists $\triangle x \in \mathcal{C}(\bar{x}) \setminus \{0\}$ such that (37) holds. Thus, we know from Lemma 22 that there exists $\bar{d} \in \mathcal{C}(\bar{x})$ such that (38) holds. Hence, for any $\tau > 0$, we have

$$\langle (\triangle x + \tau \bar{d}), \nabla^2_{xx} L(\bar{x}; \bar{y})(\triangle x + \tau \bar{d}) \rangle + \Upsilon(G'(\bar{x})(\triangle x + \tau \bar{d}))$$
$$= \langle \triangle x, \nabla^2_{xx} L(\bar{x}; \bar{y}) \triangle x \rangle + \Upsilon(G'(\bar{x})\triangle x) + \tau^2 \left( \langle \bar{d}, \nabla^2_{xx} L(\bar{x}; \bar{y}) \bar{d} \rangle + \Upsilon(G'(\bar{x})\bar{d}) \right)$$
$$+ 2\tau \left( \langle \nabla^2_{xx} L(\bar{x}; \bar{y}) \triangle x, \bar{d} \rangle + \frac{1}{2} \langle G'(\bar{x})^* \nabla \Upsilon(G'(\bar{x})\triangle x), \bar{d} \rangle \right).$$

Since $\mathcal{C}(\bar{x})$ is a convex cone, it follows from (38) that for $\tau > 0$ sufficiently small, $\triangle x + \tau \bar{d} \in \mathcal{C}(\bar{x})$ and
$$\langle (\triangle x + \tau \bar{d}), \nabla^2_{xx} L(\bar{x}; \bar{y})(\triangle x + \tau \bar{d}) \rangle + \Upsilon(G'(\bar{x})(\triangle x + \tau \bar{d})) < 0,$$

which is in contradiction with the second-order necessary condition (41). This completes the proof. $\square$

By combining Propositions 21 and 23, Lemma 18 and 22 and Theorem 17, we obtain the following result on the characterization of the (robust) isolated calmness of $S_{\text{KKT}}$.

THEOREM 24. *Let $\bar{x}$ be a feasible solution to problem (1) with $(a, b) = (0, 0)$. Suppose that the RCQ holds at $\bar{x}$. Let $\bar{y} \in M(\bar{x}, 0, 0) \neq \emptyset$. Then the following statements are equivalent:*

(i) *the SRCQ (11) holds at $\bar{x}$ with respect to $\bar{y}$ and the SOSC (14) holds at $\bar{x}$ for problem (1) with $(a, b) = (0, 0)$;*

(ii) *$\bar{x}$ is a locally optimal solution to problem (1) with $(a, b) = (0, 0)$ and $S_{\text{KKT}}$ is robustly isolated calm at the origin for $(\bar{x}, \bar{y})$;*

(iii) *$\bar{x}$ is a locally optimal solution to problem (1) with $(a, b) = (0, 0)$ and $S_{\text{KKT}}$ is isolated calm at the origin for $(\bar{x}, \bar{y})$;*

(iv) *$\bar{x}$ is a locally optimal solution to problem (1) with $(a, b) = (0, 0)$ and $F^{-1}$ is isolated calm at the origin for $(\bar{x}, \bar{y})$.*

REMARK 5. *The isolated calmness of the mapping $F^{-1}$ at the origin for $(\bar{x}, \bar{y})$ implies the following error bound result: there exist a constant $\kappa > 0$ and a neighborhood $\mathcal{V}$ of $(\bar{x}, \bar{y})$ in $\mathcal{X} \times \mathcal{Y}$ such that*

$$\|(x, y) - (\bar{x}, \bar{y})\| \leq \kappa \|F(x, y)\| \quad \forall (x, y) \in \mathcal{V}.$$

Below, we provide a nontrivial SDP example, which is modified from the example proposed by Zhou and So [39] for miminizing an unconstrained convex optimization problem penalized by the nuclear norm function, to illustrate that under the RCQ, the SOSC and the uniqueness of the Lagrange multipliers alone do not guarantee the calmness property of $S_{\text{KKT}}$ at the origin for $(\bar{x}, \bar{y})$.

EXAMPLE 3. *Let $B = \begin{bmatrix} 3/2 & -2 \\ -2 & 3 \end{bmatrix}$ and $b = B^{-1/2} \begin{bmatrix} 5/2 \\ -1 \end{bmatrix}$. Consider the following convex quadratic SDP problem:*

$$
\begin{aligned}
(P) \quad \min \quad & \tfrac{1}{2}\|x + b\|^2 + t \\
\text{s.t.} \quad & \mathcal{A}^* x + tE + I + \delta \in \mathcal{S}_+^2, \\
& t \geq 0
\end{aligned}
$$

*and its dual*

$$
\begin{aligned}
(D) \quad \max \quad & -\tfrac{1}{2}\|\mathcal{A}Y - b\|^2 - \langle I + \delta, Y \rangle \\
\text{s.t.} \quad & \langle E, Y \rangle \leq 1, \\
& Y \in \mathcal{S}_+^2,
\end{aligned}
$$

*with $\mathcal{A}Y = B^{1/2} \text{diag}(Y)$ for all $Y \in \mathcal{S}^2$, $\mathcal{A}^* x = \text{Diag}(B^{1/2} x)$ for all $x \in \Re^2$ and $\delta = \begin{bmatrix} -\varepsilon & 0 \\ 0 & \varepsilon \end{bmatrix}$ with $\varepsilon \geq 0$, where $E$ is the $2 \times 2$ symmetric matrix with elements are all one, $\text{diag}(Y)$ denotes the column vector consisting of all the diagonal entries of $Y$. The KKT system of the above convex quadratic SDP problem takes the following form:*

$$
\begin{cases}
\mathcal{A}Y(\delta) - b = x(\delta), \\
\langle E, Y(\delta) \rangle - 1 = s(\delta), \\
0 \leq t(\delta) \perp s(\delta) \leq 0, \\
\mathcal{S}_+^2 \ni Y(\delta) \perp (\mathcal{A}^* x(\delta) + t(\delta) E + I + \delta) \in \mathcal{S}_+^2,
\end{cases}
$$

whose solution set consists of $x(\delta) = B^{-1/2} \begin{bmatrix} -1+\varepsilon \\ -1-\varepsilon \end{bmatrix}$, $t(\delta) = 0$ and $Y(\delta) = \begin{bmatrix} 1+2\varepsilon & \xi \\ \xi & \varepsilon \end{bmatrix}$ with $|\xi| \leq \sqrt{\varepsilon + 2\varepsilon^2}$. When $\varepsilon = 0$, it is clear that $\bar{x} = B^{-1/2} \begin{bmatrix} -1 \\ -1 \end{bmatrix}$ and $\bar{t} = 0$ is the unique optimal solution to the unperturbed problem and $(\overline{Y}, \bar{s}) \in \mathcal{S}^2 \times \Re$ with $\overline{Y} = \begin{bmatrix} 1 & 0 \\ 0 & 0 \end{bmatrix}$ and $\bar{s} = 0$ is the unique Lagrange multiplier. Moreover, it is easy to verify that the SRCQ for the dual problem holds at $\overline{Y}$ with respect to its unique Lagrange multiplier $(0,0) \in \Re \times \mathcal{S}^2$. Thus, we know from [16, Proposition 5.3] that the SOSC for the primal problem holds at $(\bar{x}, \bar{t})$. Clearly, the calmness of $S_{\mathrm{KKT}}$ does not hold at the origin for $(\bar{x}, \bar{t}, \overline{Y}, \bar{s})$, since $\xi$ is in the order of $\varepsilon^{1/2}$. In fact, it can be checked easily that the SRCQ for the primal problem does not hold at $(\bar{x}, \bar{t})$ with respect to $(\overline{Y}, \bar{s})$.

Next, we provide an example to demonstrate that the conditions for ensuring the robust isolated calmness of $S_{\mathrm{KKT}}$ at the origin for a solution $(\bar{x}, \bar{y}) \in S_{\mathrm{KKT}}(0,0)$ are strictly weaker than the conditions for guaranteeing the Aubin property.

EXAMPLE 4. *Consider the following convex quadratic SDP problem:*

$$\begin{aligned} \min \quad & \tfrac{1}{2}(X_{11}-1)^2 + \tfrac{1}{2}(X_{22}-2X_{12})^2 \\ \text{s.t.} \quad & \langle E, X \rangle \leq 1, \\ & X \in \mathcal{S}_+^2. \end{aligned} \tag{42}$$

We know that $\overline{X} = \begin{bmatrix} 1 & 0 \\ 0 & 0 \end{bmatrix}$ is the unique optimal solution to problem (42) and $(\bar{s}, \overline{Y}) = (0,0) \in \Re \times \mathcal{S}^2$ is the unique corresponding Lagrange multiplier. Note that $H \in \mathrm{lin}(\mathcal{T}_{\mathcal{S}_+^2}(\overline{X}))$ if and only if $H = \begin{bmatrix} H_{11} & H_{12} \\ H_{12} & 0 \end{bmatrix}$ with $H_{11} \in \Re$ and $H_{12} \in \Re$. Thus, for any $(t, S) \in \Re \times \mathcal{S}^2$, we have

$$\begin{cases} \langle E, \widehat{X} \rangle = t, \\ \widehat{X} + H = S, \end{cases}$$

where $H = \begin{bmatrix} -t + \langle E, S \rangle & 0 \\ 0 & 0 \end{bmatrix} \in \mathrm{lin}(\mathcal{T}_{\mathcal{S}_+^2}(\overline{X}))$ and $\widehat{X} = S - H \in \mathcal{S}^2$, which implies that the constraint non-degeneracy (12) holds at $\overline{X}$. Note that the SOSC (14) of problem (42) at $\overline{X}$ takes the following form:

$$\langle d, \nabla_{XX}^2 L(\overline{X}; \bar{s}, \overline{Y}) d \rangle = d_{11}^2 + (d_{22} - 2d_{12})^2 > 0 \quad \forall d \in \mathcal{C}(\overline{X}) \setminus \{0\},$$

where $\mathcal{C}(\overline{X}) = \{d \in \mathcal{S}^2 \mid \langle E, d \rangle \leq 0, \; d_{22} \geq 0\}$. Therefore, we know that for any $d \in \mathcal{C}(\overline{X})$,

$$d_{11}^2 + (d_{22} - 2d_{12})^2 = 0 \implies d = 0,$$

which implies that the SOSC holds at $\overline{X}$. Thus, since the constraint non-degeneracy (12) implies the SRCQ (11), we know from Theorem 24 that $S_{\mathrm{KKT}}$ is robustly isolated calm at the origin for $(\overline{X}, \bar{s}, \overline{Y})$. Note that for any $d \in \mathcal{S}^2$ with $d_{11} = 0$ and $2d_{12} = d_{22} \neq 0$,

$$\langle d, \nabla_{XX}^2 L(\overline{X}; \bar{s}, \overline{Y}) d \rangle = d_{11}^2 + (d_{22} - 2d_{12})^2 = 0.$$

Since the affine space $\mathrm{aff}\{\mathcal{C}(\overline{X})\}$ of $\mathcal{C}(\overline{X})$ satisfies $\mathrm{aff}\{\mathcal{C}(\overline{X})\} = \mathcal{S}^2$, we know that the strong second order sufficient condition defined in [36, Definition 3.2] does not hold at $\overline{X}$. This implies the KKT point $(\overline{X}, \bar{s}, \overline{Y})$ is not a strongly regular solution of the KKT system (cf. [36, Theorem 4.1]).

Note that since the quadratic SDP problem (42) is convex, we know from [9, Proposition 5.1] (see also [24, Proposition 3.12]) that $(\overline{X}, \bar{s}, \overline{Y})$ is a strongly regular solution of the KKT system if and only if the set-valued mapping $S_{\text{KKT}}$ has the Aubin property at the origin for $(\overline{X}, \bar{s}, \overline{Y})$. Thus, this example also demonstrates that the robust isolated calmness of $S_{\text{KKT}}$ at the origin for $(\overline{X}, \bar{s}, \overline{Y})$ is a property that is strictly weaker than the Aubin property.

It is clear from the definition that the strong regularity of the KKT solution mapping $S_{\text{KKT}}$ implies both the robust isolated calmness and the Aubin property of $S_{\text{KKT}}$. Furthermore, by using Proposition 20 and Theorem 24, we can provide the following results relating to the strong regularity, the Aubin property and the robust isolated calmness of $S_{\text{KKT}}$.

COROLLARY 25. *Suppose that $\bar{x}$ is a locally optimal solution to problem* (1) *with $(a, b) = (0, 0)$ and the RCQ holds at $\bar{x}$. Let $\bar{y} \in M(\bar{x}, 0, 0) \neq \emptyset$ and $\bar{z} = G(\bar{x}) + \bar{y}$. Consider the following statements:*

(i) *The KKT point $(\bar{x}, \bar{y})$ is a strongly regular solution to the KKT system (4) with $(a, b) = (0, 0)$;*
(ii) *The mapping $\Psi^{-1}$ has the Aubin property at the origin for $(\bar{x}, \bar{z})$;*
(iii) *The KKT solution mapping $S_{\text{KKT}}$ has the Aubin property at the origin for $(\bar{x}, \bar{y})$;*
(iv) *the constraint non-degeneracy* (12) *holds at $\bar{x}$ and the SOSC* (14) *holds at $\bar{x}$ for problem* (1) *with $(a, b) = (0, 0)$;*
(v) *the SRCQ* (11) *holds at $\bar{x}$ with respect to $\bar{y}$ and the SOSC* (14) *holds at $\bar{x}$ for problem* (1) *with $(a, b) = (0, 0)$;*
(vi) *The KKT solution mapping $S_{\text{KKT}}$ is robustly isolated calm at the origin for $(\bar{x}, \bar{y})$;*
(vii) *The mapping $F^{-1}$ is isolated calm at the origin for $(\bar{x}, \bar{y})$.*

*Then it holds that*

$$(i) \implies (ii) \iff (iii) \implies (iv) \implies (v) \iff (vi) \iff (vii).$$

As mentioned in the introduction, when the non-polyhedral set $\mathcal{K}$ is the SOC or the PSD cone, under the RCQ at a locally optimal solution $\bar{x}$, the KKT solution mapping is strongly regular if and only if the constraint non-degeneracy and the strong SOSC hold [5, 36]. This, together with Example 4 and Corollary 25, implies that in order for the KKT solution mapping $S_{\text{KKT}}$ to have the Aubin property at the origin for $(\bar{x}, \bar{y})$ one must at least have both the constraint non-degeneracy (12) and a second order optimality condition that is stronger than the SOSC (14) but weaker (not necessarily strictly weaker) than the strong SOSC as being studied in [5, 36] for the nonlinear SOC programming and SDP problems[3].

**5. Conclusions.** In this paper, for a large class of canonically perturbed conic programming problems, we establish a complete characterization of the robust isolated calmness of the KKT solution mapping $S_{\text{KKT}}$, namely, under the RCQ at a locally optimal solution, the KKT solution mapping is robustly isolated calm if and only if both the SRCQ and the SOSC hold. This extension from the polyhedral sets to the case of non-polyhedral but $C^2$-cone reducible convex sets is quite satisfactory. Nevertheless, there are many other unresolved stability issues in conic optimization even if $\mathcal{K}$ is $C^2$-cone reducible. For example, the relationship between the strong regularity and the Aubin property of $S_{\mathrm{KKT}}$ is one important stability issue for the general non-convex optimization problems with non-polyhedral cone constraints. For problem (1) with a $C^2$-cone reducible convex set, at a locally optimal solution both the constraint non-degeneracy and the SOSC must hold in order to have the Aubin property. This means, at a locally optimal solution, if for the non-convex optimization problem (1) with $\mathcal{K}$ being the SOC cone or the PSD cone or something similar, the Aubin property holds for the corresponding KKT solution mapping but not the strong regularity, then a second order optimality condition weaker than the strong SOSC but stronger than the SOSC must hold. It would be interesting to know that if one can apply Mordukhovich's excellent co-derivative criterion on the Aubin property of $S_{\mathrm{KKT}}$ [34, Theorem 9.40] to derive an exact form of the second order optimality condition to be sought, together with the constraint non-degeneracy, to characterize the non-singularity of the co-derivative of Robinson's normal map $\Psi$ for understanding the Aubin property for the non-convex problem (1) with a non-polyhedral set. Also, another important unresolved issue is to provide conditions ensuring the calmness property of $S_{\mathrm{KKT}}$ for problem (1) without the uniqueness of the KKT solutions.

---
[3]In [29], the authors proved the equivalence between the Aubin property and the strong regularity for nonlinear SOC programming problem at a locally optimal solution. There was a gap in the original proof of [29, Theorem 21] on the last case (the 6th case). However, fortunately, this gap has recently been fixed in [28].

**Acknowledgments.** The authors would like to thank the two anonymous referees and the associate editor for their many valuable comments and constructive suggestions, which have substantially helped improve the quality and the presentation of this paper. In particular, the authors are in debt to one referee for suggestions on proving Proposition 9.